\documentclass[10pt]{article}	% [letterpaper, preprint, paper,11pt]{AAS}

\usepackage{amsmath,amsfonts,graphics,graphicx,bm,enumerate,epsfig,psfrag,subfig,longtable}
\usepackage{palatino,cite}
\usepackage{color,soul}
\usepackage{xcolor}
\usepackage[mathscr]{euscript} % adding \mathscr{}
\usepackage{mathabx}

\usepackage{float}
\usepackage[colorlinks=true, pdfstartview=FitV, linkcolor=black, citecolor= black, urlcolor= black]{hyperref}
\usepackage{footnpag}			       % make footnote symbols restart on each page

\usepackage{framed,fancybox}
\definecolor{shadecolor}{gray}{0.99}
\usepackage{multirow}
\usepackage{booktabs}
\usepackage{colortbl}
\definecolor{Gray1}{gray}{0.85}

\newenvironment{shadedframe}{%
 \MakeFramed {\FrameRestore}}
{\endMakeFramed}

\newcommand{\m}[1]{{\bf{#1}}}

\newcommand{\bb}[1]{\mathbb #1}

\newcommand{\C}[1]{{\cal {#1}}}

\newcommand{\singlespacing}{\renewcommand{\baselinestretch}{1}\normalsize\normalfont}

\date{}

\begin{document}

\title{\bf  Mesh Refinement Method for Solving Optimal \\ Control Problems with Nonsmooth Solutions \\ Using Jump Function Approximations}

\author{Alexander T.~Miller\thanks{Ph.D.~Student, NDSEG Fellow, Department of Mechanical and Aerospace Engineering.  E-mail:  alexandertmiller@ufl.edu.} \\ William W. Hager\thanks{Distinguished Professor, Department of Mathematics. E-mail: hager@ufl.edu} \\ Anil~V.~Rao\thanks{Professor, Erich Farber Faculty Fellow, and University Term Professor, Department of Mechanical and Aerospace Engineering.  E-mail:  anilvrao@ufl.edu.  Associate Fellow AIAA.  Corresponding Author. \newline \newline It is noted that a preliminary version of this research was presented at the 2018 AIAA Guidance, Navigation, and Control Conference.} \\\\ {\em University of  Florida} \\ {\em Gainesville, FL 32611}}

\maketitle{} 		
\singlespacing

%-------------------------
\begin{abstract}
A mesh refinement method is described for solving optimal control problems using Legendre-Gauss-Radau collocation.  The method detects discontinuities in the control solution by employing an edge detection scheme based on jump function approximations.  When discontinuities are identified, the mesh is refined with a targeted $h$-refinement approach whereby the discontinuity locations are bracketed with mesh points.  The remaining smooth portions of the mesh are refined using previously developed techniques.  The method is demonstrated on two examples, and results indicate that the method solves optimal control problems with discontinuous control solutions using fewer mesh refinement iterations and less computation time when compared with previously developed methods.
\end{abstract}

%--------------------------
\section{Introduction}
Over the past few decades, direct collocation methods have become increasingly popular for solving optimal control problems numerically.  Direct collocation methods are state and control parameterization methods where the dynamics are approximated at a set of specially chosen points called {\em collocation points}.  The optimal control problem is then transcribed to a finite-dimensional nonlinear programming problem (NLP) \cite{Betts3,Hargraves1}.  The NLP is then is solved numerically using well known software \cite{Gill1,Biegler2}.  In any direct collocation method, the optimal control problem is approximated on a {\em mesh}, where the mesh is a division of the independent variable into segments called {\em mesh intervals} over which collocation is performed.  Traditional direct collocation methods take the form of an $h$ method where a low-order method is employed (typically, a method such as trapezoidal or Runge-Kutta is used) and the order of the method is the same in every mesh interval.   Accuracy in an $h$ method is then achieved by increasing the number of mesh intervals and/or adjusting the locations of the intervals \cite{Betts3,Jain1,Betts10}.  More recently, research has explored $p$ methods.  In a $p$ method, the order of the method is varied in each mesh interval, but the number of mesh intervals remains small.  Accuracy using a $p$ method is then achieved by increasing the order of the approximation in each mesh interval.  In order to achieve maximum effectiveness, $p$ methods have been developed using {\em Gaussian quadrature} collocation \cite{Benson2,Garg1,Elnagar1}.  Gauss quadrature collocation methods employ Legendre-Gauss\cite{Benson2} (LG), Legendre-Gauss-Radau\cite{Garg1} (LGR), or Legendre-Gauss-Lobatto\cite{Elnagar1} (LGL) points and converge at an exponential rate when the solution is smooth and well-behaved \cite{HagerHouRao15a,HagerHouRao15b,HagerHouRao15c,HagerHouRao16a,HagerMohapatraRao16b}.

Although $h$ methods have been used extensively and $p$ methods have shown promise on problems where the solution can be approximated accurately using a low-degree or moderate-degree polynomial, both the $h$ and $p$ approaches have limitations.  In the case of an $h$ method, obtaining a high-accuracy solution typically requires the use of an extremely fine mesh.  In the case of a $p$ method, obtaining a high-accuracy solution typically requires the use of an unreasonably large-degree polynomial.  To overcome the limitations with $h$ and $p$ methods, recent research has focused on $hp$ methods.  As the name suggests, $hp$ methods adjust both the number and placement of mesh intervals and the order of the approximation in each mesh interval. While $hp$ methods were originally developed as finite-element methods for solving partial differential  equations \cite{Babuska3,Babuska4,Babuska5,Babuska1,Babuska2}, in the past decade a framework has been established for solving optimal control problems numerically using $hp$ methods \cite{Darby2,Darby3,Patterson2015,Liu2015,Peng1,Liu2018}.  Recent research has shown that $hp$ direct collocation methods can outperform $h$ or $p$ methods both in terms of computational efficiency and in reducing the size of the finite-dimensional approximation.

When solving an optimal control problem computationally using either an $h$, $p$, or $hp$ method, mesh refinement is employed to improve the accuracy of the discretization \cite{Betts10}.  In optimal control, mesh refinement produces a sequence of meshes on which the problem is discretized and the corresponding NLP is solved.  The solution on a given mesh is generally higher in accuracy than the solution on the preceding mesh, and mesh refinement terminates when a specified accuracy tolerance has been attained.  The decision-making process employed to iteratively refine the mesh is known as a {\em mesh refinement method}.

Various mesh refinement methods for solving optimal control problems have been developed previously.  Reference~\citen{Betts3} develops a primarily $h$ method which subdivides the mesh based on an error estimate for the state.  A $ph$ strategy is described in Ref.~\citen{Patterson2015} where priority is placed on increasing the polynomial degree first, but mesh intervals are subdivided and set to a low degree approximation if the required polynomial degree exceeds a specified maximum.  Reference~\citen{Liu2018} introduces a method which analyzes the decay rates of Legendre polynomial coefficients to determine whether mesh intervals are smooth or nonsmooth, and then takes $p$ or $h$ refinement actions accordingly.  Reference~\citen{Gong3} describes a method that uses a differentiation matrix to identify and place knots at locations where the control rate of change is high.  Reference~\citen{Zhao2} employs density functions to allocate grid points such that the integrated density is equal across each mesh interval.  Beyond the references just listed, the reader is referred to\cite{Cuthrell1,Cuthrell2,Schlegel2005,Jain1,Darby2,Darby3,Liu2015,Peng1,Zhao2019,Munos,Grune,Rannacher2009,Rannacher2012,Dorao1,Dorao2,Heinrichs1,Galvao1,Karniadakis1,Mitchell2014} and the references therein for further information regarding mesh refinement methods for optimal control and $hp$ refinement techniques more generally.

Although many of the aforementioned mesh refinement methods can be effective when the solution is smooth, these methods can be less effective when solving an optimal control problem with a nonsmooth solution.  When the solution is nonsmooth, the piecewise-smooth nature of the solution may not match the piecewise-smooth parameterization of the solution.  Any loss in solution accuracy due to such a mismatch must be overcome by generating an appropriately partitioned mesh.  Thus, a key challenge in mesh refinement is to detect and accurately locate the nonsmooth features of the solution and then to adjust the mesh appropriately.  The new mesh must adequately capture the piecewise-smooth nature of the solution while simultaneously maintaining computational efficiency by employing the sparsest mesh possible.

It is desirable that a mesh refinement method exactly identify all locations of nonsmoothness in the solution and then partition the mesh at those locations.  In such a method, the mesh would contain a mesh point that coincides with each location of nonsmoothness.  The transcription of the optimal control problem on such a mesh would result in a piecewise-smooth approximation that matches the piecewise-smooth behavior of the solution.  Current mesh refinement methods for optimal control, however, make little attempt to identify locations of nonsmoothness or do so with low accuracy.  In addition, typical mesh refinement strategies for refining the mesh once a nonsmooth feature is identified tend to result in repeated subdivision of the mesh around the nonsmooth location until the accuracy tolerance is attained.  Such mesh refinement approaches can be computationally intensive because it may take many mesh refinement iterations to achieve the desired accuracy.  In addition, the mesh resulting from such an approach may place an unnecessarily large number of mesh intervals near the location of nonsmoothness, producing a larger discretization than necessary.

Motivated by the prevalence of nonsmooth behavior in the solutions to optimal control problems and the need for better mesh refinement methods to handle such solutions, this paper describes a new method for solving optimal control problems whose solutions are nonsmooth.  In the method of this paper, discontinuities in the control solution are identified and located using jump function approximations \cite{Gelb1999, Yoon2005, Gelb2006, Gelb2007, Yoon2008}.  Jump function approximations are an effective tool in detecting discontinuities, because they tend towards zero everywhere except at the discontinuity locations where they tend towards the value of the jump in the underlying function.  While most research on jump function approximations has focused on their mathematical foundations and with image and signal processing in mind (see Ref.~\citen{Gelb1999, Yoon2005, Gelb2006, Gelb2007, Yoon2008}), this paper extends the use of jump function approximations to optimal control.  In particular, this paper employs jump function approximations for the purposes of discontinuity detection in a novel mesh refinement method.  Whenever discontinuities are identified, the mesh is refined with a targeted $h$-refinement approach which brackets each of the identified discontinuity locations with mesh points.  On subsequent mesh refinement iterations, the mesh points that bracket a discontinuity are reused and their locations updated to reflect the higher accuracy and precision of the discontinuity location estimate on the latest mesh.  Such an approach confines discontinuities to small mesh intervals rapidly, and does so in a way which does not add unnecessary size to the mesh (and therefore the NLP).

The contributions of this paper are as follows.  The use of jump function approximations is extended to optimal control.  Specifically, jump function approximations are incorporated into the mesh refinement process as a tool for detecting and accurately locating discontinuities using only the numerical solution on a given mesh.  The next contribution is the mesh refinement method itself.  The method combines an effective discontinuity detection scheme with a specialized $h$-refinement procedure that rapidly increases the resolution of the mesh around identified discontinuities while keeping the mesh size relatively small.  Computational benefits gained by employing the new method are demonstrated on two examples as a final contribution.

The remainder of this paper is organized as follows.  Section~\ref{sect:Bolza} introduces the optimal control problem written in Bolza form.  A brief review of jump function approximations is given in Section~\ref{sect:jfBackground}.  Section~\ref{sect:LGR} describes the transcription of the optimal control problem to an NLP.  This is followed by Sections~\ref{sect:meshRefinementMethod} and \ref{sect:algorithm} which describe the mesh refinement method and the corresponding mesh refinement algorithm respectively.  The method is then demonstrated on two examples in Section~\ref{sect:examples}.  Finally, Sections~\ref{sect:discussion}~and~\ref{sect:conclusion} discuss the results and draw conclusions about the method.

% ----------------------------------------------------------------
\section{Bolza Optimal Control Problem \label{sect:Bolza}}
Consider the following optimal control problem written in Bolza form.  Determine the initial and final times, $t_0$ and $t_f$, as well as the state, $\m{y}(\tau)\in\bb{R}^{n_y}$, and the control, $\m{u}(\tau)\in\bb{R}^{n_u}$, on the domain $\tau \in [-1,+1]$ that minimize the cost functional
\begin{equation}\label{eq:bolza-cost-tau}
  \C{J} = \C{M}(\m{y}(-1),t_0,\m{y}(+1),t_f)+\frac{t_f-t_0}{2}\int_{-1}^{+1}
\C{L}(\m{y}(\tau),\m{u}(\tau), t(\tau, t_0, t_f))\, d\tau,
\end{equation}
while satisfying the state dynamics relations
\begin{equation}\label{eq:bolza-dyn-tau}
  \frac{d\m{y}}{d\tau} - \frac{t_f-t_0}{2}\m{a}(\m{y}(\tau),\m{u}(\tau), t(\tau, t_0, t_f) ) = \m{0}, 
\end{equation}
the boundary conditions
\begin{equation}\label{eq:bolza-bc-tau}
  \m{b}(\m{y}(-1), t_0, \m{y}(+1), t_f ) \leq \m{0},
\end{equation}
and the path constraints
\begin{equation}\label{eq:bolza-path-tau}
  \m{c}(\m{y}(\tau),\m{u}(\tau), t(\tau, t_0, t_f) ) \leq \m{0},
\end{equation}
where the functions $\C{M}$, $\C{L}$, $\m{a}$, $\m{b}$, and $\m{c}$ are defined by the mappings
\begin{equation*} \label{eq:bolza-mappings}
  \begin{array}{rl}
	\C{M} &: \bb{R}^{n_y} \times \bb{R} \times \bb{R}^{n_y} \times \bb{R} \to \bb{R}, \\
	\C{L} &:~\bb{R}^{n_y} \times \bb{R}^{n_u} \times \bb{R} \to \bb{R}, \\
	\m{a} &:~\bb{R}^{n_y} \times \bb{R}^{n_u} \times \bb{R} \to \bb{R}^{n_y}, \\
	\m{b} &:~\bb{R}^{n_y} \times \bb{R} \times \bb{R}^{n_y} \times \bb{R} \to \bb{R}^{n_b}, \\
	\m{c} &:~\bb{R}^{n_y} \times \bb{R}^{n_u} \times \bb{R} \to \bb{R}^{n_c}, \\
  \end{array}
\end{equation*}
and the affine relation 
\begin{equation}\label{tau-to-t}
  t \equiv t(\tau,t_0,t_f) = \frac{t_f-t_0}{2}\tau + \frac{t_f+t_0}{2},
\end{equation}
relates the computational domain $\tau \in [-1,+1]$ to the time interval $t \in [t_0,t_f]$.

% -----------------------------------------------------------------
\section{Jump Function Approximations\label{sect:jfBackground}}
Solutions to the Bolza optimal control problem defined in Section~\ref{sect:Bolza} are often nonsmooth.  Nonsmoothness in the solution may take the form of a discontinuous control (for example, a bang-bang control), a discontinuous state or control derivative, or discontinuities in higher-order derivatives.  Nonsmooth behavior in the solution can be difficult to approximate numerically, because numerical methods for optimal control typically assume a smooth or piecewise smooth parameterization which often does not match the piecewise smooth nature of the solution.  The mismatch between the solution and the parameterization of the solution results in a decrease in the numerical solution accuracy in the neighborhood of the discontinuity.  Therefore, in order to improve the accuracy of a numerical solution to an optimal control problem with a nonsmooth solution, the locations of the associated discontinuities must first be accurately determined so that the mesh can be refined appropriately.

The method of this paper employs jump function approximations for the purposes of estimating discontinuities in the solution of an optimal control problem.  A brief background on jump functions and methods to approximate jump functions is provided in Section~\ref{subsect:jfBackground}.  With the appropriate background established, Section~\ref{subsect:JFmotivation4OCP} discusses how jump function approximations are useful when applied to solving optimal control problems whose solutions are nonsmooth.

\subsection{Background\label{subsect:jfBackground}}
A jump function is defined as follows.  Let $f: \bb{R} \to \bb{R}$ be an arbitrary function defined on the interval $t \in [t_0,t_f]$.  The jump function of $f(t)$, denoted $[f](t)$, is defined as
\begin{equation} \label{eq:def-JF}
	[f](t) = f(t^+) - f(t^-),
\end{equation}
where $f(t^+)$ and $f(t^-)$ are the right-hand and left-hand limits of $f(t),~t \in (t_0,t_f)$.  Equation~\eqref{eq:def-JF} dictates that the jump function be zero across intervals where $f(t)$ is continuous and that the jump function take on the value of the jump in $f(t)$ at those locations where $f(t)$ is discontinuous.  Such a property is useful in discontinuity detection, because one can discern the discontinuity locations by observing where the jump function is non-zero.  However, constructing the jump function via Eq.~\eqref{eq:def-JF} requires the value of the underlying function, $f(t)$, be known over the entire domain of interest.

Jump function approximations offer a more practical approach for determining discontinuity locations.  Much of the research in developing jump function approximations has been conducted in the context of digital image processing where discontinuities take the form of edges in a two-dimensional image.  Here, a brief overview of the major developments in jump function approximation techniques is provided, but is restricted to one dimension.  Further information regarding jump function approximations and edge detection can be found in Ref.~\citen{Gelb1999, Yoon2005, Gelb2006, Gelb2007, Yoon2008} and the references therein.

The original jump function approximation is derived from a Fourier analysis of the underlying function as follows.  Consider an arbitrary function $f: \bb{R} \to \bb{R}$ with domain $t \in [-\pi,\pi]$.  The truncated Fourier series expansion of $f(t)$ is given by
\begin{equation} \label{eq:FourierSeries}
	S_N [f](t) = \frac{a_0}{2} + \sum_{n=1}^{N} \left[a_n \cos \left(n t\right) + b_n \sin \left(n t\right) \right],
\end{equation}
where $a_n$ and $b_n$ are the Fourier coefficients.  The conjugate Fourier series expansion of $f(t)$ is then expressed as
%with Fourier coefficients defined by
%\begin{equation} \label{eq:FourierCoef}
%  \begin{array}{rclcl}
%	a_n &=& \frac{1}{\pi} \int_{-\pi}^{\pi} f(t) \cos (n t) dt & , & n = 0,\ldots,N,\\
%	b_n &=& \frac{1}{\pi} \int_{-\pi}^{\pi} f(t) \sin (n t) dt & , & n = 1,\ldots,N.
%  \end{array}
%\end{equation}
%The conjugate Fourier series expansion of $f(t)$ is then given by
\begin{equation} \label{eq:ConFourierSeries}
	\tilde{S}_N [f](t) = \sum_{n=1}^{N} \left[a_n \sin \left(n t\right) - b_n \cos \left(n t\right) \right].
\end{equation}
For real-valued, continuous functions the truncated Fourier series approaches $f(t)$ and the truncated conjugate Fourier series approaches zero as the number of terms increases.  In contrast, real-valued functions that contain discontinuities produce a truncated Fourier series that experiences Gibbs phenomenon at the discontinuity locations.  Likewise, the truncated conjugate Fourier series does not approach zero at the discontinuity locations in $f(t)$ even as $N$ approaches infinity.  Instead, work by Luk{\'a}cs \cite{Bary1,Zygmund1} has shown that the conjugate Fourier series satisfies the condition
\begin{equation} \label{eq:LukacsResult}
	\lim_{N \to \infty} \left\lbrace \frac{-\pi}{\log N} \tilde{S}_N [f](t) \right\rbrace = [f](t),
\end{equation}
where $[f](t)$ is the jump function.  In other words, the truncated conjugate Fourier series scaled by negative $\pi / \log{N}$ {\em approximates} the jump function.  It is noted that the result of Eq.~\eqref{eq:LukacsResult} has been extended to include {\em generalized conjugate Fourier partial sums} of the form
\begin{equation} \label{eq:GenConFourierSeries}
	\tilde{S}^{\sigma}_N [f](t) = \sum_{n=1}^{N} \sigma \left(\frac{n}{N} \right) \left[a_n \sin \left(n t\right) - b_n \cos \left(n t\right) \right],
\end{equation}
where $\sigma_{n,N} = \sigma \left(\frac{n}{N}\right)$ are known as {\em concentration factors} \cite{Gelb1999}.  For admissible concentration factors, a generalized conjugate Fourier partial sum has the property
\begin{equation} \label{eq:GelbResult}
	\lim_{N \to \infty} \tilde{S}^{\sigma}_N [f](t) = [f](t).
\end{equation}
The discrete analog to Eq.~\eqref{eq:GenConFourierSeries} is constructed from a uniform grid of samples $f(t_j)$, with grid points defined by $t_j = -\pi + (j + N)\Delta t$ and $\Delta t = 2 \pi / (2N + 1)$.  These samples of $f(t)$ permit the formation of a {\em generalized discrete conjugate Fourier partial sum} given by
\begin{equation} \label{eq:GenConFourierSeries-discrete}
	\tilde{T}^{\gamma}_N [f](t) = \sum_{n=1}^{N} \gamma \left(\frac{n}{N} \right) \left[A_n \sin \left(n t\right) - B_n \cos \left(n t\right) \right],
\end{equation}
%where $\{A_n,B_n\}$ are the discrete Fourier coefficients
% defined by
%\begin{equation} \label{eq:FourierCoef-discrete}
%  \begin{array}{rclcl}
%	A_n &=& \frac{\Delta t}{\pi} \sum\limits_{j = -N}^{N} f(t_j) \cos (n t_j) & , & n = 0,\ldots,N,\\
%	B_n &=& \frac{\Delta t}{\pi} \sum\limits_{j = -N}^{N} f(t_j) \sin (n t_j) & , & n = 1,\ldots,N,
%  \end{array}
%\end{equation}
where $\{A_n,B_n\}$ are the discrete Fourier coefficients and $\gamma_{n,N} = \gamma \left(\frac{n}{N} \right)$ are the discrete concentration factors \cite{Gelb1999}.  Like its continuous counterpart, the generalized discrete conjugate Fourier partial sum shares the property of Eq.~\eqref{eq:GelbResult} and converges to the jump function as the number of terms approaches infinity.

An entirely different jump function approximation which does not rely upon a Fourier analysis of the underlying function was introduced in Ref.~\citen{Yoon2005}.  The approximation is developed as follows.  Consider an arbitrary function $f(t)$ defined on $t \in [t_0,t_f]$ and assume the function is sampled on a grid of $N$ points, not necessarily evenly spaced.  Let $\C{S}_{t} = \{ t_1,\ldots,t_{m+1} \}$ be the set of $m+1$ closest grid points surrounding $t$, and let $\C{S}_{t}^{+}$ contain only the elements of $\C{S}_{t}$ which are larger than $t$.  The jump function of $f(t)$ is approximated by
\begin{equation} \label{eq:Lmf}
	L_m f(t) = \frac{1}{q_m (t)} \sum_{t_j \in \C{S}_t} c_j(t) f(t_j) \approx [f](t),
\end{equation}
where $q_m (t)$ is defined by
\begin{equation} \label{eq:qm}
	q_m (t) = \sum_{t_j \in \C{S}_{t}^{+}} c_j(t),
\end{equation}
$c_j (t)$ is defined by
\begin{equation} \label{eq:cj}
	c_j (t) = \frac{m!}{\prod\limits_{\substack{i = 1 \\i \neq j}}^{m + 1} (t_j - t_i)},
\end{equation}
and $m$ specifies the order of the approximation.  Higher order approximations converge to the jump function faster outside the neighborhood of discontinuities but have oscillatory behavior in the vicinity of discontinuities.  The oscillations are reduced via the \textbf{minmod} function, defined here as
\begin{equation} \label{eq:MMLmf}
MM\left( L_{\mathscr{M}} f(t) \right) = \left\{
  \begin{array}{ll}
    \min\limits_{m \in \mathscr{M}} L_m f(t) & if~L_m f(t) > 0~~\forall~m \in \mathscr{M}, \\
    \max\limits_{m \in \mathscr{M}} L_m f(t) & if~L_m f(t) < 0~~\forall~m \in \mathscr{M}, \\
    0                           & otherwise,\\
  \end{array}\right.
\end{equation}
where $\mathscr{M} \subset \bb{N}^{+}$ is a finite set of choices of the approximation order $m$.  Equation~\eqref{eq:MMLmf} converges to the jump function as the grid point spacing in $\C{S}_{t}$ approaches zero \cite{Yoon2005}.

%Convergence of Eq.~\eqref{eq:MMLmf} to the jump function is derived in Ref.~\citen{Yoon2005} as
%\begin{equation} \label{eq:MMLmf-convergence}
%MM\left( L_{\mathscr{M}} f(t) \right) = \left\{
%  \begin{array}{ll}
%    [f](\zeta) + \C{O}(h(t)) & if~t_{j-1} \leq \zeta,t \leq t_j, \\
%    \C{O}(h^{\min (\mathscr{M}_t, k)}(t)) & if~f \in C^k(I_t),
%  \end{array}\right.
%\end{equation}
%where
%\begin{equation}
%h(t) = \max \{ |t_i - t_{i-1}|~:~t_{i-1},t_i \in \C{S}_t \},
%\end{equation}
%\begin{equation}
%\mathscr{M}_t = \max \{ m \in \mathscr{M}~|~\# \C{S}_t = m+1,~I_t \cap J = \emptyset \},
%\end{equation}
%$I_t$ is the smallest closed interval such that $\C{S}_t \subset I_t$ with $\# \C{S}_t \leq \mathscr{M}_t + 1$, and $\zeta \in J$ where $J$ is the set of all discontinuities in $f(t)$ on $[t_0,t_f]$.  As can be seen, $MM\left( L_{\mathscr{M}} f(t) \right)$ approximates $[f](t)$ for finite grids and converges to $[f](t)$ as the grid point spacing approaches zero.

The convergence properties of Equations~\eqref{eq:GenConFourierSeries},~\eqref{eq:GenConFourierSeries-discrete}, and~\eqref{eq:MMLmf} to the jump function allow one to distinguish between the continuous regions of $f(t)$ and intervals of $f(t)$ where a discontinuity is suspected.  Neighborhoods surrounding discontinuities in $f(t)$ will see jump function approximations which approach the value of the jump.  Conversely, the jump function approximation values will tend towards zero outside the neighborhoods of discontinuities.  Therefore, one can detect discontinuities by observing where the absolute value of an appropriate jump function approximation exceeds a specified threshold.  In this manner, jumps in $f(t)$ of sufficiently large magnitude cause the jump function approximation to exceed the threshold in the neighborhood of the jump and that neighborhood is identified as being likely to contain a discontinuity.
%It may be tempting to choose a low threshold so that smaller jumps may be detected.  However, choosing a threshold which is too small can result in continuous behavior (such as regions of rapid change) to be misidentified as jump behavior.

% Backup in case you want to scrap the background section and simply summarize with references:
%In the case of Reference~\citen{Gelb1999}, $[f](t)$ is estimated via a \textit{generalized conjugate Fourier partial sum} which utilizes uniformly-spaced samples of the underlying function $f(t)$.  Another approach, described in Reference~\citen{Yoon2005}, approximates jump functions using only local stencils which can be non-uniformly spaced.  Yet other approaches exist, but for the purpose of brevity, we refer the reader to References~\citen{Gelb1999, Yoon2005, Gelb2006, Gelb2007} and the references therein for further information.

\subsection{Application to Optimal Control\label{subsect:JFmotivation4OCP}}
Jump function approximations are an effective tool for estimating the locations of nonsmooth behavior in the solutions of optimal control problems.  In particular, a jump function approximation built from the optimal control problem solution allows one to make an assessment as to whether the underlying function is continuous or discontinuous and where the discontinuities (if any) might lie.  In principle, jump function approximations could be produced to detect a variety of nonsmooth behavior such as control discontinuities, control derivative discontinuities, or corners in the state.  The question of which jump function to approximate and which method to approximate it with is dependent upon the type of nonsmooth  behavior one seeks to identify and the (perhaps limited) information available about the underlying function in the numerical solution.

Typically, the nonsmooth behavior of interest occurs in the state and/or the control solution.  It is generally the case that the state has more continuous derivatives than the control, and that nonsmooth features in the state solution tend to coincide in location with nonsmooth features of the control solution.  Therefore, searching for nonsmoothness in the control solution is usually sufficient and easier to do.  Of the types of nonsmooth control solution behavior, discontinuous controls tend to be the biggest cause for concern when solving an optimal control problem numerically.  For these reasons and for the purposes of clarity and brevity, the scope of this paper is restricted to searching for discontinuities in the control.

The strategy for determining control discontinuities requires first generating jump function approximations for each of the control components.  While any one of the methods discussed in Section~\ref{subsect:jfBackground} could be applied to generate jump function approximations for the control components, all three [Eqs.~\eqref{eq:GenConFourierSeries}, \eqref{eq:GenConFourierSeries-discrete}, and \eqref{eq:MMLmf}] require a particular set of information.  In the case of Eq.~\eqref{eq:GenConFourierSeries}, the Fourier coefficients used in the jump function approximation require a continuous representation of the control.  In contrast, Eq.~\eqref{eq:GenConFourierSeries-discrete} requires the control be sampled on an evenly spaced grid in order to obtain the discrete Fourier coefficients.  Finally, Eq.~\eqref{eq:MMLmf} only requires local samples of the control which do not need to be evenly spaced.

In this research, the numerical solution to the optimal control problem is attained by employing a direct collocation method which produces values for the control at the collocation points.  The spacing of the collocation points across $[t_0,t_f]$ is method dependent, but uneven spacing is typical.  Therefore, the jump function approximation of Eq.~\eqref{eq:MMLmf} is most practical from an implementation standpoint, because Eq.~\eqref{eq:GenConFourierSeries-discrete} would require interpolation to an evenly-spaced grid, and Eq.~\eqref{eq:GenConFourierSeries} would require the discrete control values be converted to some continuous parameterization.  Thus, the remainder of this paper employs Eq.~\eqref{eq:MMLmf} exclusively for approximating the jump functions of the control components.

%------------------------------------------------------------
\section{Legendre-Gauss-Radau Collocation\label{sect:LGR}}
Section~\ref{sect:jfBackground} introduced three jump function approximation techniques and discussed how Eq.~\eqref{eq:MMLmf} is particularly useful in direct collocation methods for generating jump function approximations of the control components.  In principle, Eq.~\eqref{eq:MMLmf} could be applied within a number of direct collocation schemes.  Here, we restrict ourselves to Legendre-Gauss-Radau (LGR) collocation.  The remainder of this paper employs LGR collocation exclusively to illustrate how jump function approximations can be incorporated into the existing framework of direct collocation methods for optimal control and the computational benefits of doing so.

A number of factors weigh in on the decision to use an LGR collocation scheme.  The LGR method allows for a highly general problem formulation, takes advantage of the integration accuracy and exponential convergence rates obtained by employing a Gaussian quadrature method, and can be posed as an integration method which allows for a convenient approach to estimate the state error.  Additionally, a convergence theory has been established \cite{HagerHouRao15c,HagerHouRao16a}.  The accuracy of the LGR method is particularly important, because the accuracy of each control component's jump function approximation is inherently tied to the accuracy of the numerical solution for the control.

\subsection{Multiple Interval Bolza Formulation \label{subsect:bolza-segmented}}
In the $hp$ discretization of the LGR collocation method, the domain of the Bolza problem described in Section~\ref{sect:Bolza}, $\tau\in[-1,+1]$, is partitioned into a {\em mesh} consisting of $K$ {\em mesh intervals}.  The mesh intervals are defined as $\C{S}_k=[T_{k-1},T_k],\; k=1,\ldots,K$, where $-1 = T_0 < T_1 < \ldots < T_K = +1$.  Together, the mesh intervals satisfy the property that $\displaystyle \bigcup_{k=1}^{K} \C{S}_k=[-1,+1]$.  Let $\m{y}^{(k)}(\tau)$ and $\m{u}^{(k)}(\tau)$ be the state and control in $\C{S}_k$. The Bolza optimal control problem of Eqs.~\eqref{eq:bolza-cost-tau}-\eqref{eq:bolza-bc-tau} is expressed in multiple interval form as follows.  Minimize the cost functional 
\begin{equation}\label{eq:bolza-cost-segmented}
    \C{J}  = \C{M}(\m{y}^{(1)}(-1),t_0,\m{y}^{(K)}(+1),t_f)
     + \frac{t_f-t_0}{2}\sum_{k=1}^K \int_{T_{k-1}}^{T_k}
    \C{L}(\m{y}^{(k)}(\tau),\m{u}^{(k)}(\tau),t (\tau, t_0, t_f) )\, d\tau,
\end{equation}
subject to the state dynamics relations
\begin{equation}\label{eq:bolza-dyn-segmented}
    \frac{d\m{y}^{(k)}(\tau)}{d\tau} - \frac{t_f-t_0}{2}\m{a}(\m{y}^{(k)}(\tau),\m{u}^{(k)}(\tau),t (\tau, t_0, t_f) ) = \m{0},   \quad (k=1,\ldots,K),
\end{equation}
the boundary conditions
\begin{equation} \label{eq:bolza-bc-segmented}
\m{b}(\m{y}^{(1)}(-1),t_0,\m{y}^{(K)}(+1),t_f) \leq  \m{0},  
\end{equation}
the path constraints
\begin{equation} \label{eq:bolza-path-segmented}
    \m{c}(\m{y}^{(k)}(\tau),\m{u}^{(k)}(\tau),t (\tau, t_0, t_f) ) \leq \m{0},    \quad (k=1,\ldots,K),
\end{equation}
and the state continuity constraints,
\begin{equation} \label{eq:bolza-segmented-stateContinuity}
\m{y}(T_k^{+}) - \m{y}(T_k^{-}) = \m{0},    \quad (k=1,\ldots,K-1).
\end{equation}

\subsection{Formation of the Nonlinear Program \label{subsect:NLP}}
The segmented Bolza problem of Eqs.~\eqref{eq:bolza-cost-segmented}--\eqref{eq:bolza-segmented-stateContinuity} is discretized using collocation at LGR points \cite{Garg1,Garg2,Garg3,Kameswaran1,Patterson2015}.  The state is approximated on each mesh interval by a Lagrange polynomial with support points at the Legendre-Gauss-Radau (LGR) nodes \cite{Abramowitz1}, $\left(\tau_1^{(k)},\ldots,\tau_{P_k}^{(k)}\right) \in [T_{k-1},T_k)$, and the non-collocated endpoint, $\tau_{P_k+1}^{(k)} = T_k$.  The resulting $P_k$-degree polynomial state parameterization is expressed on each mesh interval as
\begin{equation}\label{eq:LGR-stateApprox}
\m{y}^{(k)}(\tau)  \approx \m{Y}^{(k)}(\tau) = \sum_{j=1}^{P_k+1} \m{Y}_{j}^{(k)}
\ell_{j}^{(k)}(\tau), \quad  \ell_{j}^{(k)}(\tau) = \prod_{\stackrel{l=1}{l\neq j}}^{P_k+1}\frac{\tau-\tau_{l}^{(k)}}{\tau_{j}^{(k)}-\tau_{l}^{(k)}},    \quad (k=1,\ldots,K),
\end{equation}
where $\tau \in [-1,+1]$, and $\ell_{j}^{(k)}(\tau),$ $j=1,\ldots,P_k+1$ are a basis of Lagrange polynomials.  Differentiating Eq.~\eqref{eq:LGR-stateApprox} with respect to $\tau$ leads to
\begin{equation}\label{eq:LGR-stateDerivApprox}
  \frac{d\m{Y}^{(k)}(\tau)}{d\tau} = \sum_{j=1}^{P_k+1}\m{Y}_{j}^{(k)}\frac{d\ell_j^{(k)}(\tau)}{d\tau},    \quad (k=1,\ldots,K).
\end{equation}
The state derivative approximation of Eq.~\eqref{eq:LGR-stateDerivApprox} is collocated with the right-hand side of the system dynamics at the LGR points of each mesh interval, producing the state dynamics approximation,
\begin{equation} \label{eq:LGR-dynamics}
 \sum_{j=1}^{P_k+1}D_{ij}^{(k)} \m{Y}_j^{(k)} = \frac{t_f-t_0}{2}\m{a}(\m{Y}_i^{(k)},\m{U}_i^{(k)},t (\tau_i^{(k)},t_0, t_f)),   \quad (i=1,\ldots,P_k, \quad k=1,\ldots,K),
\end{equation}
where
\begin{displaymath}
  D_{ij}^{(k)} = \frac{d\ell_j^{(k)}(\tau_i^{(k)})}{d\tau}, \quad
   (i=1,\ldots,P_k, \quad j=1,\ldots,P_k+1),
\end{displaymath}
are the elements of the $P_k\times (P_k+1)$ {\em Legendre-Gauss-Radau differentiation matrix} \cite{Garg1} in mesh interval $\C{S}_k$.  The remaining constraints are similarly discretized, forming the discrete boundary conditions,
\begin{equation} \label{eq:LGR-bc} 
   \m{b}(\m{Y}_{1}^{(1)},t_0,\m{Y}_{P_K+1}^{(K)},t_f)  \leq \m{0}, 
\end{equation}
and the discretized path constraints,
\begin{equation}\label{eq:LGR-path}
   \m{c}(\m{Y}_{i}^{(k)},\m{U}_{i}^{(k)},t (\tau_i^{(k)},t_0, t_f)) \leq \m{0},    \quad (i=1,\ldots,P_k, \quad k=1,\ldots,K).
\end{equation}
Finally, the cost functional is approximated as a cost function by applying LGR quadrature rules to estimate the integral portion of Eq.~\eqref{eq:bolza-cost-segmented} on each mesh interval.  The cost function is defined as
\begin{equation}\label{eq:LGR-cost}
    \C{J}  \approx \C{M}(\m{Y}_{1}^{(1)},t_0,\m{Y}_{P_K+1}^{(K)},t_f)
     +   \frac{t_f-t_0}{2}  \sum_{k=1}^{K} \sum_{j=1}^{P_k} 
  w_{j}^{(k)} \C{L}(\m{Y}_{j}^{(k)},\m{U}_{j}^{(k)},t (\tau_j^{(k)}, t_0, t_f)),
\end{equation}
where $w_{j}^{(k)}$ is the $j^{th}$ LGR weight in mesh interval $\C{S}_k$.

The discretization of Eqs.~\eqref{eq:bolza-dyn-segmented}--\eqref{eq:bolza-path-segmented} and the quadrature approximation of Eq.~\eqref{eq:bolza-cost-segmented} combine to form the following nonlinear program.  Minimize the cost function of Eq.~\eqref{eq:LGR-cost} subject to the discrete dynamics constraints of Eq.~\eqref{eq:LGR-dynamics}, the boundary conditions of Eq.~\eqref{eq:LGR-bc}, and the discretized path constraints of Eq.~\eqref{eq:LGR-path}.  Notice that the state continuity constraints of Eq.~\eqref{eq:bolza-segmented-stateContinuity} are implicitly applied by using the same variable for $Y_{P_k + 1}^{(k)}$ and $Y_{1}^{(k+1)}$ at each of the interior mesh points.  Lastly, it is noted that
\begin{equation}\label{eq:LGR-numColPts}
  P = \sum_{k=1}^{K} P_k,
\end{equation}
is the total number of LGR points.

\subsection{Approximation of Solution Error \label{subsect:solError}}
Suppose the nonlinear program of Eqs.~\eqref{eq:LGR-dynamics}--\eqref{eq:LGR-cost} has been solved on a mesh, $\C{S}_k,~k = 1,\ldots,K$, with $P_k$ collocation points in mesh interval $\C{S}_k$.  An estimate of the discretization error on the current mesh must be obtained in order to assess the accuracy of the solution.  The approach used in this paper estimates the relative error in the state solution as a proxy for measuring the discretization error.  The method employed is identical to Ref.~\citen{Patterson2015} and is summarized here.

The objective is to approximate the error in the state at a set of $M_k = P_k + 1$ LGR points $\left(\hat{\tau}_1^{(k)},\ldots,\hat{\tau}_{M_k}^{(k)}\right)$, where 
$\hat{\tau}_1^{(k)}=\tau_1^{(k)}=T_{k-1}$, and $\hat{\tau}_{M_k+1}^{(k)}=T_k$.  Let the values of the state approximation at the points $\left(\hat{\tau}_1^{(k)},\ldots,\hat{\tau}_{M_k}^{(k)}\right)$ be
denoted $\left(\m{y}(\hat{\tau}_1^{(k)}),\ldots,\m{y}(\hat{\tau}_{M_k}^{(k)})\right)$.  Next, let the control be approximated in mesh interval $\C{S}_k$ with the Lagrange polynomial
\begin{equation}\label{eq:control-approximation-RPM}
\m{U}^{(k)}(\tau) = \sum_{j=1}^{P_k} \m{U}_{j}^{(k)}
  \hat{\ell}_{j}^{(k)}(\tau), \quad  \hat{\ell}_{j}^{(k)}(\tau) = \prod_{\stackrel{l=1}{l\neq j}}^{P_k}\frac{\tau-\tau_{l}^{(k)}}{\tau_{j}^{(k)}-\tau_{l}^{(k)}},
\end{equation}
and let the control approximation at $\hat{\tau}_i^{(k)}$ be denoted $\m{u}(\hat{\tau}_i^{(k)})$ for $i = 1,\ldots, M_k$.  The value of the right-hand side of the dynamics at $(\m{Y}(\hat{\tau}_i^{(k)}), \m{U}(\hat{\tau}_i^{(k)}), \hat{\tau}_i^{(k)})$ is used to construct an improved approximation of the state.  Let $\hat{\m{Y}}^{(k)}$ be a polynomial of degree at most $M_k$ that is defined on the interval ${\cal S}_k$. If the derivative of $\hat{\m{Y}}^{(k)}$ matches the dynamics at each of the Radau quadrature points $\hat{\tau}_i^{(k)}$, $i = 1,\ldots, M_k$, then we have 
\begin{equation}\label{eq:mesh-Radau-integral-approximation}
  \hat{\m{Y}}^{(k)}(\hat{\tau}_{j+1}^{(k)}) =
\m{Y}^{(k)}(\hat{\tau}_1^{(k)})+\frac{t_f-t_0}{2}
\sum_{l=1}^{M_k}\hat{I}_{jl}^{(k)}\m{a}
\left(\m{Y}^{(k)}(\hat{\tau}_l^{(k)}),\m{U}^{(k)}(\hat{\tau}_l^{(k)}),t(\hat{\tau}_l^{(k)}, t_0, t_f)\right), \quad\begin{array}{c}j=1,\ldots,M_k,\end{array}
\end{equation}
where $\hat{I}_{jl}^{(k)},\; j,l=1,\ldots,M_k$, is the $M_k\times M_k$ LGR integration matrix corresponding to the LGR points defined by $\left(\hat{\tau}_1^{(k)},\ldots,\hat{\tau}_{M_k}^{(k)}\right)$.  Comparing the interpolated values, $\m{Y}(\hat{\tau}_l^{(k)}),~l=1,\ldots,M_k+1$, with the integrated values, $\hat{\m{Y}}(\hat{\tau}_l^{(k)}),~l=1,\ldots,M_k+1$, the {\em absolute} and {\em relative errors} in the $i^{th}$ component of the state at $(\hat{\tau}_1^{(k)},\ldots,\hat{\tau}_{M_k+1}^{(k)})$ are defined, respectively, as
\begin{equation}\label{eq:absolute-and-relative-errors}
  \begin{array}{lcl}
    E_i^{(k)}(\hat{\tau}_l^{(k)})&=&\left|\hat{Y}_i^{(k)}(\hat{\tau}_l^{(k)})-Y_i^{(k)}(\hat{\tau}_l^{(k)})\right|,\vspace{6pt}\\
   e^{(k)}_i(\hat{\tau}_l^{(k)}) & = &\displaystyle
    \frac{E^{(k)}_i(\hat{\tau}_l^{(k)})}{\displaystyle1+\max_{\stackrel{j \in[1,\ldots,P_k+1]}{k \in[1,\ldots,K]}}\left|Y_i^{(k)}(\tau_j^{(k)})\right|},
    \end{array} \quad \left[\begin{array}{c}l=1,\ldots,M_k+1,\\i=1,\ldots,n_y,\end{array}\right].
\end{equation}
The {\em maximum relative error} in mesh interval $\C{S}_k$ is then defined as
\begin{equation}\label{eq:maximum-relative-error}
 e^{(k)}_{\max} =  \max_{\stackrel{i\in[1,\ldots,n_y]}{l\in[1,\ldots,M_k+1]}}e^{(k)}_i(\hat{\tau}_l^{(k)}).
\end{equation}

%-----------------------------------------------------------
\section{$hp$-Adaptive Mesh Refinement Method\label{sect:meshRefinementMethod}}
Suppose the nonlinear program of Eqs.~\eqref{eq:LGR-dynamics}--\eqref{eq:LGR-cost} has been solved on a mesh $\C{S}_k = [T_{k-1},T_k],~k = 1,\ldots,K$, with $P_k$ LGR points in mesh interval $\C{S}_k$.  Suppose further that the maximum relative error estimate, $e^{(k)}_{\max}$ of Eq.~\eqref{eq:maximum-relative-error}, exceeds a desired error tolerance, $\epsilon$, on one or more of the mesh intervals.  In order to satisfy the desired error tolerance on each mesh interval, the current mesh must be refined using an appropriate mesh refinement method.  The LGR method is then applied to the new mesh and the cycle repeated until the error tolerance is satisfied.

This section develops a novel $hp$-adaptive mesh refinement method which detects and brackets any control discontinuities detected in the solution.  The method refines the mesh in a two pass sequence referred to as "nonsmooth mesh refinement" and "smooth mesh refinement" respectively.  The novelty of the method lies in the nonsmooth mesh refinement method and its ability to interface with existing mesh refinement methods which are used to carry out smooth mesh refinement.  The idea is to detect, locate, and bracket nonsmooth features of the solution (control discontinuities in this research) during nonsmooth mesh refinement, and then refine the remaining smooth portions of the mesh with a chosen smooth mesh refinement method.  In this manner, existing mesh refinement methods which perform poorly when discontinuities are present in the solution can be combined with the nonsmooth mesh refinement approach to further improve their performance.

The nonsmooth and smooth refinement methods are described next in Section~\ref{subsect:nonsmooth-refine} and Section~\ref{subsect:smooth-refine} respectively.  Throughout this process, it is useful to define a few terms.  The mesh used to discretize and solve the optimal control problem on the current iteration is called the \textit{current mesh} while the mesh obtained after the current iteration of mesh refinement is complete is called the \textit{new mesh}.  The new mesh is used to discretize and solve the optimal control problem again, becoming the current mesh on the next iteration.  The term \textit{intermediate mesh} describes the mesh as it transitions from the current mesh to the new mesh.

%At each stage of the nonsmooth refinement process an index mapping variable $s_{k}$ is used to relate the ${k}^{th}$ mesh interval on the intermediate mesh back to the $s_{{k}}^{th}$ mesh interval on the current mesh so that smooth refinement can be carried out later on using the solution and error estimates obtained on the current mesh.

Each stage of the mesh refinement process is guided by whether and where the solution is thought to be smooth or nonsmooth.  The terms \textit{smooth segment} and \textit{nonsmooth segment} are introduced here to label particular regions of the mesh as follows.  A nonsmooth segment is a region of the mesh defined by $[T_m,T_{m+2}],~m \in \{0,\ldots,K-2\}$ which bounds one (and only one) of the identified discontinuities.  The remaining portions of the mesh are labeled as smooth segments, because they contain none of the identified discontinuities.  Note that each nonsmooth segment is comprised of the two mesh intervals which form the bracket around the identified discontinuity, and each smooth segment spans the mesh intervals which connect one nonsmooth segment to the next.  Moreover, the initial mesh is comprised of a single smooth segment encompassing the entire mesh on $[T_0,T_K]$, and the smooth and nonsmooth segment labels are updated during the nonsmooth mesh refinement process.

\subsection{Nonsmooth Mesh Refinement\label{subsect:nonsmooth-refine}}
Nonsmooth mesh refinement detects, locates, and then brackets discontinuities identified on the current mesh.  In this research, only control discontinuities are considered and their locations are estimated using jump function approximations.  Nonsmooth mesh refinement begins by employing the method of Section~\ref{subsubsect:disc-detect} to identify and estimate the locations of the control discontinuities.  After discontinuity locations have been estimated, a sequence of specialized refinement actions are taken to bracket newly identified discontinuities, update the brackets of previously identified discontinuities, and relabel any nonsmooth segments which are now thought to be smooth.

The nonsmooth mesh refinement actions occur in the following sequence.  Suppose the method of Section~\ref{subsubsect:disc-detect} detects a total of $n_d$ discontinuities whose estimated locations are $d_1 < d_2 < \ldots < d_{n_d}$ with estimated uncertainty bounds $\{d_i^-,d_i^+\},~i = 1,\ldots,n_d$.  These uncertainty bounds are first adjusted using the method of Section~\ref{subsubsect:normalize} in order to prevent overlap with one another and to limit the estimated uncertainty bounds to a single smooth or nonsmooth segment.  Next, any newly identified discontinuities located on a smooth segment of the current mesh are bracketed via the method of Section~\ref{subsubsect:disc-bracket}.  The newly identified discontinuities may instead lie on a nonsmooth segment of the current mesh, indicating a previously identified discontinuity has been identified again.  In such a case, the existing bracket is updated by employing the method of Section~\ref{subsubsect:bracket-update}.  Finally, the method of Section~\ref{subsubsect:bracket-relabel} relabels any nonsmooth segments which are no longer thought to contain a discontinuity.  Detailed descriptions of the discontinuity detection procedure and the nonsmooth mesh refinement actions are described next.

\subsubsection{Discontinuity Detection\label{subsubsect:disc-detect}}
Discontinuities are detected on a case-by-case basis for each control component using jump function approximations of the form given in Eq.~\eqref{eq:MMLmf} constructed from the numerical control solution.  Consider the solution for the $i^{th}$ control component $U_{i}(\tau_j^{(k)}),~i \in \{1,\ldots,n_u\},~j=1,\ldots,P_k,~k=1,\ldots,K$.  The range of the control component solution is first normalized to $[0,1)$ via the transformation
\begin{equation}
u_i(\tau_j^{(k)}) = \frac{\displaystyle U_i(\tau_j^{(k)}) - U_{i,\min}}
{\displaystyle 1 + U_{i,\max} - U_{i,\min}},
\end{equation}
where $U_{i,\min}$ and $U_{i,\max}$ are the minimum and maximum values of the control component solution.  The normalized control component solution, $u_i(\tau_j^{(k)})$, and the corresponding collocation points, $\tau_j^{(k)}$ (on $[-1,1)$ of the entire mesh), are applied to Eqs.~\eqref{eq:Lmf}--\eqref{eq:MMLmf}, producing a jump function approximation for the normalized control component solution.  The jump function approximation is only evaluated at $\tau_{j+\frac{1}{2}}^{(k)} = \frac{1}{2}\left(\tau_j^{(k)} + \tau_{j+1}^{(k)}\right),$ $j = 1,\ldots,P_k$ and $k \in \{1,\ldots,K\}$ such that the error estimate, $e_{\max}^{(k)}$ of Eq.~\ref{eq:maximum-relative-error}, exceeds the specified error tolerance ($e_{\max}^{(k)} > \epsilon$).  In other words, the jump function approximation is only evaluated on those mesh intervals which need to be refined.  Now, let the evaluation of the jump function approximation at $\tau_{j+\frac{1}{2}}^{(k)}$ be denoted by $MM\left( L_{\mathscr{M}}~u_i(\tau_{j+\frac{1}{2}}^{(k)}) \right)$.  A discontinuity is considered present at $\tau_{j+\frac{1}{2}}^{(k)}$ when the condition
\begin{equation}\label{eq:discontinuityThreshold}
\left| MM\left( L_{\mathscr{M}}~u_i(\tau_{j+\frac{1}{2}}^{(k)}) \right) \right| \geq \eta,
\end{equation}
is satisfied for at least one of the normalized control components ($u_i,~i \in \{1,\ldots,n_u\}$).  Note that the parameter $0 < \eta < 1$ is a user-defined threshold which specifies the relative size of jumps that should be detected.  Smaller values of $\eta$ are more likely to identify discontinuities present in the control solution, but run a higher risk of attaining false positives.  Likewise, larger values of $\eta$ reduce the likelihood of false positives, but run a higher risk of attaining false negatives.

Bounds on the uncertainties in the discontinuity locations are also estimated.  Suppose Eq.~\eqref{eq:discontinuityThreshold} is satisfied at $\tau_{j+\frac{1}{2}}^{(k)}$ for some $j \in \{1,\ldots,P_k\}$ and $k \in \{1,\ldots,K\}$.  This indicates that a discontinuity is present in the normalized control solution somewhere on $\tau \in [\tau_{j}^{(k)},\tau_{j+1}^{(k)}]$ \cite{Yoon2005}.  However, the goal is to identify discontinuities in the optimal control, not its numerical approximation.  Due to the additional uncertainty incurred by employing the numerical control solution as a proxy for the optimal control, a safety factor, $\mu \geq 1$, is introduced to extend the estimated discontinuity uncertainty bounds.  The more conservative bounds estimates are defined by
\begin{equation}\label{eq:disc-boundsEstimates}
\begin{array}{rcl}
d^- & = & \tau_{j+\frac{1}{2}}^{(k)} - \mu \left( \tau_{j+\frac{1}{2}}^{(k)} - \tau_{j}^{(k)} \right), \\
d^+ & = & \tau_{j+\frac{1}{2}}^{(k)} + \mu \left( \tau_{j+1}^{(k)} - \tau_{j+\frac{1}{2}}^{(k)} \right).
\end{array}
\end{equation}
Larger values of $\mu$ are more likely to produce bounds which contain the discontinuity in the optimal control.  However, smaller values of $\mu$ increase the resolution around the discontinuity when the problem is re-solved on the new mesh, as will become clear later.

\subsubsection{Adjust Discontinuity Uncertainty Bounds\label{subsubsect:normalize}}
Suppose the smooth or nonsmooth segment $[T_m,T_n],~0 \leq m < n \leq K$ on the current mesh contains the newly identified discontinuities $d_i,~i=1,\ldots,D$ where $1 \leq D \leq n_d$ and $T_m < d_1 < \ldots < d_D < T_n$.  The associated uncertainty bounds $\{d_i^-,d_i^+\},~i=1,\ldots,D$ are adjusted as follows.  The exterior uncertainty bounds, $d_1^-$ and $d_D^+$, are first confined to $[T_m,T_n]$ by redefining them as
\begin{equation}
  \begin{array}{rcl}
    d_1^- &=& \max(T_m,d_1^-),\\
    d_D^+ &=& \min(T_n,d_D^+).
  \end{array}
\end{equation}
Likewise, any overlapping interior uncertainty bounds ($d_i^+ > d_{i+1}^-,~i \in \{ 1,\ldots,D-1 \}$) are redefined as
\begin{equation}
    d_i^+ = d_{i+1}^- = \frac{1}{2}\left(d_i + d_{i+1}\right),
\end{equation}
in order to resolve overlap.

\subsubsection{Bracket Discontinuities\label{subsubsect:disc-bracket}}
\begin{figure}[h]
  \centering
  \includegraphics[height= 2 in]{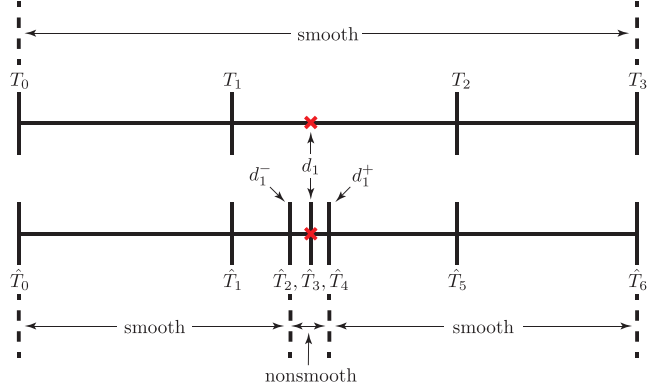}
  \caption{Example of the bracketing technique for a single discontinuity (red x).  The intermediate mesh is shown before (top) and after (bottom) refinement.\label{fig:bracket}}
\end{figure}

Suppose the intermediate mesh is currently comprised of mesh intervals $\C{S}_k = [T_{k-1},T_k],~k = 1,\ldots,K$, with $P_k$ LGR points in mesh interval $\C{S}_k$.  Suppose further that the smooth segment $[T_m,T_n],~0 \leq m < n \leq K$ contains the newly identified discontinuities $d_i,~i=1,\ldots,D$, with estimated uncertainty bounds $\{d_i^-,d_i^+\},~i=1,\ldots,D$, where $1 \leq D \leq n_d$ and $T_m < d_1 < \ldots < d_D < T_n$.  The intermediate mesh is refined as follows.  For each discontinuity, $d_i,~i=1,\ldots,D$, the intermediate mesh on $[d_i^-,d_i^+]$ is replaced by a new nonsmooth segment consisting of two mesh intervals.  The new mesh intervals span $[d_i^-,d_i]$ and $[d_i,d_i^+]$ respectively and are allocated four collocation points each.  All remaining mesh intervals are either unaffected by the change or are subdivided in the process.  In either case, these mesh intervals retain the same number of collocation points as the mesh interval they originated from before the bracketing operation.  Figure~\ref{fig:bracket} demonstrates the process for a single discontinuity.

% Nonsmooth Refinement Example: bracket discontinuities

\subsubsection{Update Brackets\label{subsubsect:bracket-update}}
Suppose the intermediate mesh is currently comprised of mesh intervals $\C{S}_k = [T_{k-1},T_k],~k = 1,\ldots,K$, with $P_k$ LGR points in mesh interval $\C{S}_k$.  Suppose further that the nonsmooth segment $[T_m,T_{m+2}],~0 \leq m < m+2 \leq K$, contains the newly identified discontinuities $d_i,~i=1,\ldots,D$, with associated uncertainty bounds $\{d_i^-,d_i^+\},~i=1,\ldots,D$, where $1 \leq D \leq n_d$ and $T_m < d_1 < \ldots < d_D < T_{m+2}$.  The expected scenario is $D = 1$, because each nonsmooth segment is constructed to bracket a single discontinuity.  However, $D > 1$ is possible when the locations of discontinuities lie close together relative to the resolution provided by the current mesh on the previous iteration of mesh refinement.  The general case where $D \geq 1$ is treated here for completeness.

The discontinuity bracket is updated as follows.  New brackets are formed by replacing the nonsmooth segment on $[T_m,T_{m+2}]$ by $D$ nonsmooth segments spanning $[d_i^-,d_i^+],~i = 1,\ldots,D$, respectively.  Each new nonsmooth segment consists of two mesh intervals covering $[d_i^-,d_i]$ and $[d_i,d_i^+]$ respectively and are allocated four collocation points each.  Each connecting interval $[d_i^+,d_{i+1}^-],~i \in \{ 1,\ldots,D-1 \}$, such that $d_i^+ < d_{i+1}^-$ forms a new mesh interval which is allocated four collocation points and is labeled as a smooth segment.  Finally, the mesh intervals $\C{S}_m = [T_{m-1},T_m]$ and $\C{S}_{m+3} = [T_{m+2},T_{m+3}]$ are extended to include $[T_m,d_1^-]$ and $[d_D^+,T_{m+2}]$ respectively, noting that the number of allocated collocation points remains unchanged.  However, in the special case where $\C{S}_m$ or $\C{S}_{m+3}$ is contained on a nonsmooth segment or does not exist ($m = 0$ or $m+2 = K$), a new smooth segment is created instead.  The new smooth segment is comprised of a single mesh interval spanning $[T_m,d_1^-]$ or $[d_D^+,T_{m+2}]$ respectively and is allocated four collocation points.  Note that no new smooth segment is created when $T_m = d_1^-$ or $d_D^+ = T_{m+2}$ respectively.

\begin{figure}[h]
  \centering
  \includegraphics[height=2 in]{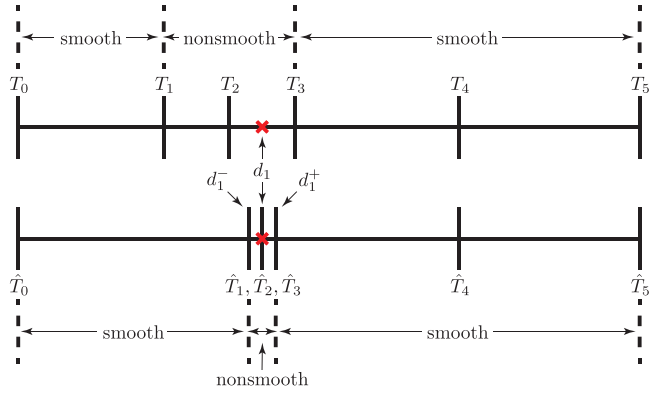}
  \caption{Example of the update technique for a single discontinuity (red x).  The intermediate mesh is shown before (top) and after (bottom) refinement.\label{fig:update}}
\end{figure}

Figure~\ref{fig:update} demonstrates the update process for a single discontinuity.  As can be seen, the net effect of the update procedure is to contract the bracket around the discontinuity location.  Two benefits are apparent.  The first is that the resolution around the discontinuity location increases with each successive bracket update.  A second benefit is that the mesh points that bracket the discontinuity are reused, thereby limiting the growth in size of the mesh.

% Nonsmooth Refinement Example: update brackets

\subsubsection{Relabel Brackets\label{subsubsect:bracket-relabel}}
\begin{figure}[h]
  \centering
  \includegraphics[height=2 in]{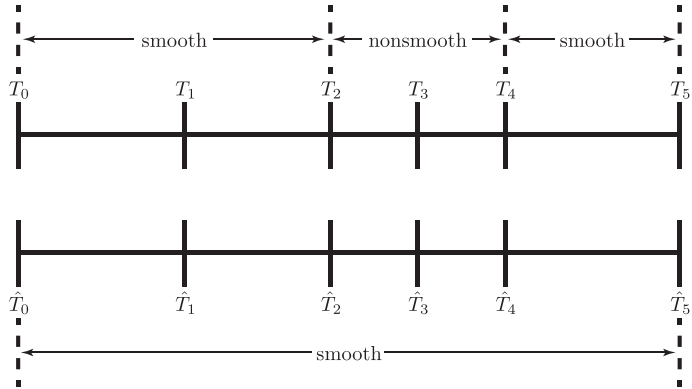}
  \caption{Example of the relabeling procedure.  The intermediate mesh is shown before (top) and after (bottom) the segments have been relabeled.\label{fig:relabel}}
\end{figure}

Suppose the intermediate mesh is currently comprised of mesh intervals $\C{S}_k = [T_{k-1},T_k],~k = 1,\ldots,K$, with $P_k$ collocation points in mesh interval $\C{S}_k$.  Suppose further that none of the newly identified discontinuities, $d_i,~i=1,\ldots,n_d$, are contained on the nonsmooth segment $[T_m,T_{m+2}],~0 \leq m < m+2 \leq K$, and that either $\C{S}_{m+1} = [T_m,T_{m+1}]$ or $\C{S}_{m+2} = [T_{m+1},T_{m+2}]$ requires refinement (that is, either $e_{\max}^{(m+1)} > \epsilon$ or $e_{\max}^{(m+2)} > \epsilon$).  The nonsmooth segment is relabeled as smooth and is combined with any adjacent smooth segments.  However, no structural changes are made to the intermediate mesh as can be seen in Fig.~\ref{fig:relabel} which illustrates the process.

% Nonsmooth Refinement Example: relabel brackets

\subsection{Smooth Mesh Refinement\label{subsect:smooth-refine}}
Suppose the current mesh is comprised of mesh intervals $\C{S}_k = [T_{k-1},T_k],~k = 1,\ldots,K$, with $P_k$ LGR points in mesh interval $\C{S}_k$.  Suppose further that the nonsmooth mesh refinement method of Section~\ref{subsect:nonsmooth-refine} has produced an intermediate mesh consisting of mesh intervals $\hat{\C{S}}_{\hat{k}} = [\hat{T}_{\hat{k}-1},\hat{T}_{\hat{k}}],~\hat{k} = 1,\ldots,\hat{K}$, with $\hat{P}_{\hat{k}}$ collocation points in mesh interval $\hat{\C{S}}_{\hat{k}}$.  In its current state, the intermediate mesh has bracketed all of the identified discontinuities, leaving only those mesh intervals contained on smooth segments of the intermediate mesh to be refined.  While many methods could fulfill the smooth mesh refinement role, this research refines each smooth segment using one of three different methods currently available in the literature\cite{Darby3,Patterson2015,Liu2018}.

It may be ambiguous at first glance how one should proceed with smooth mesh refinement given the changes made during nonsmooth mesh refinement.  The state and control solution as well as the maximum relative error estimates, $e_{\max}^{(k)},~k = 1,\ldots,K$, were all obtained using the current mesh, not the intermediate mesh.  In order to proceed, a mesh interval index mapping heuristic is employed to relate the index of any mesh interval contained on a smooth segment of the intermediate mesh back to a unique index of a mesh interval from the current mesh.  The map $\hat{k} \rightarrow k$ allows mesh refinement decisions to be made for $\hat{\C{S}}_{\hat{k}}$ based upon the solution and error estimates obtained for $\C{S}_k $.

The mesh interval indices are mapped as follows.  Only mesh intervals contained on a smooth segment of the intermediate mesh are being refined, so only those mesh intervals are considered in the map.  Three possibilities arise for the mapping $\hat{k} \rightarrow k$.  \textbf{Case 1:} $\hat{\C{S}}_{\hat{k}} \bigcap (T_{k-1},T_k) \neq \emptyset$ for a unique choice of $\hat{k}$ and $k$ such that $\C{S}_k$ lies on a smooth segment of the current mesh.  \textbf{Case 2:} $\hat{\C{S}}_{\hat{k}} \equiv \C{S}_k$ for a unique choice of $\hat{k}$ and $k$ such that $\C{S}_k$ lies on a nonsmooth segment of the current mesh.  \textbf{Case 3:} neither case 1 nor case 2 is satisfied for a particular choice of $\hat{k}$, no map is attained, and the mesh interval $\hat{\C{S}}_{\hat{k}}$ is not refined any further.  Note that the third case only arises when the update procedure of Section~\ref{subsubsect:bracket-update} produces a new smooth segment, meaning $\hat{\C{S}}_{\hat{k}}$ has already been refined during the nonsmooth refinement process.

After obtaining the map $\hat{k} \rightarrow k$, one last step is necessary to avoid over-refinement.  The error estimates $e_{\max}^{(k)},~k \in \{1,\ldots,K\}$ are set to zero when the corresponding mesh interval, $\C{S}_k$, contains one or more of the identified discontinuities, $d_i,~i \in \{1,\ldots,n_d\}$.  The reasoning is that mesh refinement has already been performed on $\C{S}_k$ during nonsmooth refinement, so additional refinement is overkill.  Following this final adjustment, smooth mesh refinement is performed one smooth segment at a time using the method of Ref.~\citen{Patterson2015}, Ref.~\citen{Liu2018}, or Ref.~\citen{Darby3} until refinement is complete and the new mesh has been formed.

%\clearpage
% -----------------------------------------------------------------
\section{Mesh Refinement Algorithm\label{sect:algorithm}}
An overview of the mesh refinement algorithm appears below.  The mesh refinement iteration is denoted by $M$ and is incremented by one with each loop of the algorithm.  The algorithm terminates when either the mesh error tolerance, $\epsilon$, is satisfied on each mesh interval or when $M$ reaches a prescribed limit, $M_{max}$.

\begin{shadedframe}
\vspace{-10pt}
\begin{center}
 \shadowbox{\bf Mesh Refinement Method}
\end{center}
\begin{enumerate}[{\bf Step 1:}]
\item Set $M=0$ and supply initial mesh.  All mesh intervals form a single, smooth segment at the start.
\item Solve Radau collocation NLP of Section~\ref{subsect:NLP} on mesh $M$. \label{step:solve}
\item Compute maximum relative error $e^{(k)}_{\max}$ in $\mathcal{S}_k,\;k=1,\ldots,K$, using Eq.~\eqref{eq:maximum-relative-error}.
\item If $e^{(k)}_{\max} \leq\epsilon$ for all $k\in \{ 1,\ldots,K \}$ or
  $M>M_{\max}$, then quit.  Otherwise, proceed to {\bf Step \ref{step:disc-detect}}.  
\item Apply the nonsmooth mesh refinement method of Section~\ref{subsect:nonsmooth-refine}.
	\begin{enumerate}[{\bf (a):}]
	\item Identify and estimate the locations of discontinuities using the method of Section~\ref{subsubsect:disc-detect}. \label{step:disc-detect}
	\item Adjust the discontinuity location uncertainty bounds using the method of Section~\ref{subsubsect:normalize}.
	\item Bracket new discontinuity locations using the method of Section~\ref{subsubsect:disc-bracket}.
	\item Update existing discontinuity brackets using the method of Section~\ref{subsubsect:bracket-update}
	\item Relabel discontinuity brackets using the method of Section~\ref{subsubsect:bracket-relabel}.
	\end{enumerate}
\item Apply the smooth mesh refinement method of Section~\ref{subsect:smooth-refine}.
\item Increment $M$ by unity and return to {\bf Step \ref{step:solve}}.  
\end{enumerate}
\end{shadedframe}

%\clearpage
%----------------------------------------------
\section{Examples\label{sect:examples}}
In this section the mesh refinement method described in Section~\ref{sect:meshRefinementMethod} is demonstrated on two examples.  The first example has a discontinuous control solution which highlights the ability of the method to detect, locate, and  bracket discontinuities quickly and efficiently.  The second example has a continuous control solution and demonstrates the ability of the method to discern when discontinuities are not present.  In evaluating the performance of the method developed in this paper, comparisons will be made against the method of Ref.~\citen{Liu2018}.  While comparisons could also be made against other mesh refinement methods (for example, the methods of Ref.~\citen{Patterson2015} and Ref.~\citen{Darby3}), the method of Ref.~\citen{Liu2018} tends to outperform these previously developed methods and offers a more competitive comparison.

When using the various methods, the terminology $hp$ and $hp$-$(\mu)$ is adopted to refer to the method of Ref.~\citen{Liu2018} without and with nonsmooth mesh refinement respectively, where $\mu \geq 1$ denotes the value of the safety factor used in Eq.~\eqref{eq:disc-boundsEstimates}.  Results are shown for three safety factors $\mu = (1,~1.5,~2)$ to illustrate the performance of the method as a function of the safety factor.  In addition, jump function approximation orders $\mathscr{M} = \{1,\ldots,6\}$ are used to evaluate Eq.~\eqref{eq:MMLmf} and the discontinuity detection threshold is set to $\eta = 0.1$.  All results were obtained using the MATLAB optimal control software $\mathbb{GPOPS-II}$ \cite{Patterson2014} running with the NLP solver IPOPT \cite{Biegler2} in full Newton mode with the linear solver MA57 \cite{Duff1}, an optimality tolerance of $10^{-9}$, and mesh refinement accuracy tolerances of $\epsilon = (10^{-6}, 10^{-7}, 10^{-8})$.  All first and second derivatives were supplied to IPOPT using the built-in sparse central differencing method in $\mathbb{GPOPS-II}$ that uses the method of Ref.~\cite{Patterson2012}.  Next, for each result shown the initial mesh consists of ten uniformly-spaced mesh intervals with four collocation points in each mesh interval.  Moreover, the initial guess consists of a straight line connecting the given (or guessed) initial and terminal values of the variables.  All computations were performed on a 2.4 GHz 8-Core Intel Core i9 MacBook Pro running macOS Catalina version 10.15.1 with 32 GB of 2400 MHz DDR4 RAM and MATLAB version R2019b.  The CPU times reported in this paper are 20-run averages of the execution time.

% Robot Arm----------------------------------------
\subsection{Example 1: Minimum Time Reorientation of a Robotic Arm\label{subsect:robotArm}}
\begin{figure}[h]
  \centering
  \begin{tabular}{cc}
  \includegraphics[height=2.25in]{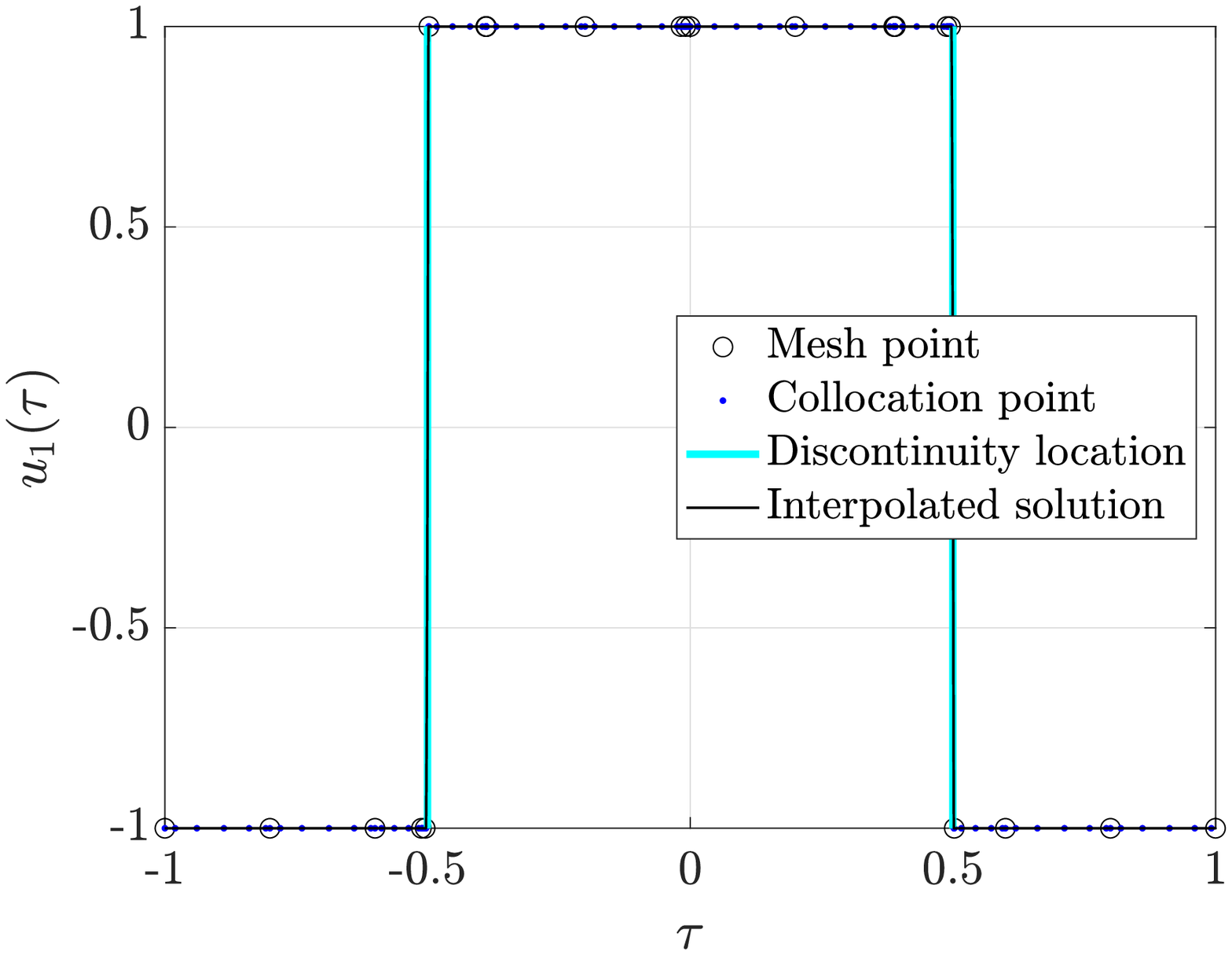}
  &
  \includegraphics[height=2.25in]{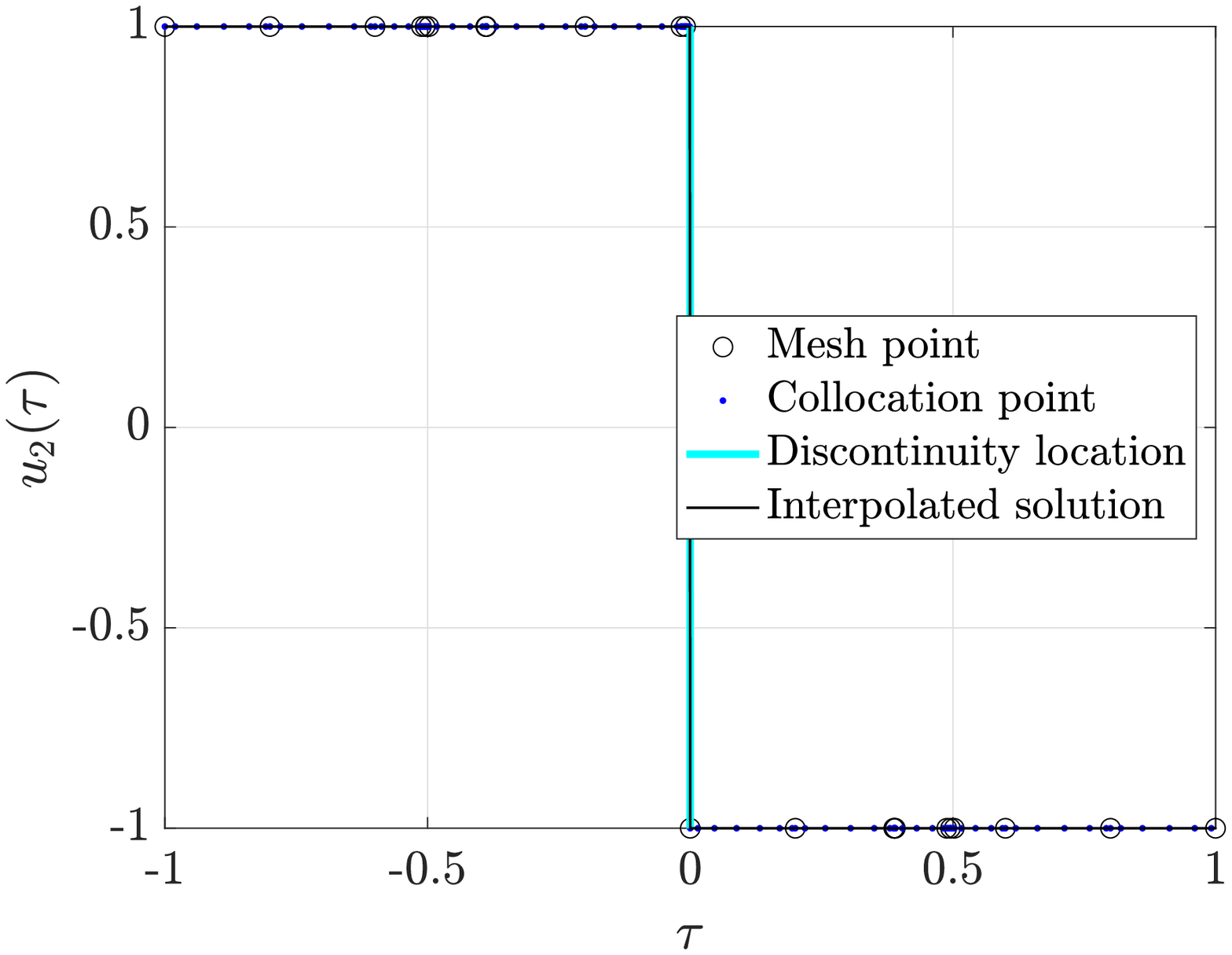}
  \end{tabular}
  \\\vspace{1mm}
  \includegraphics[height=2.25in]{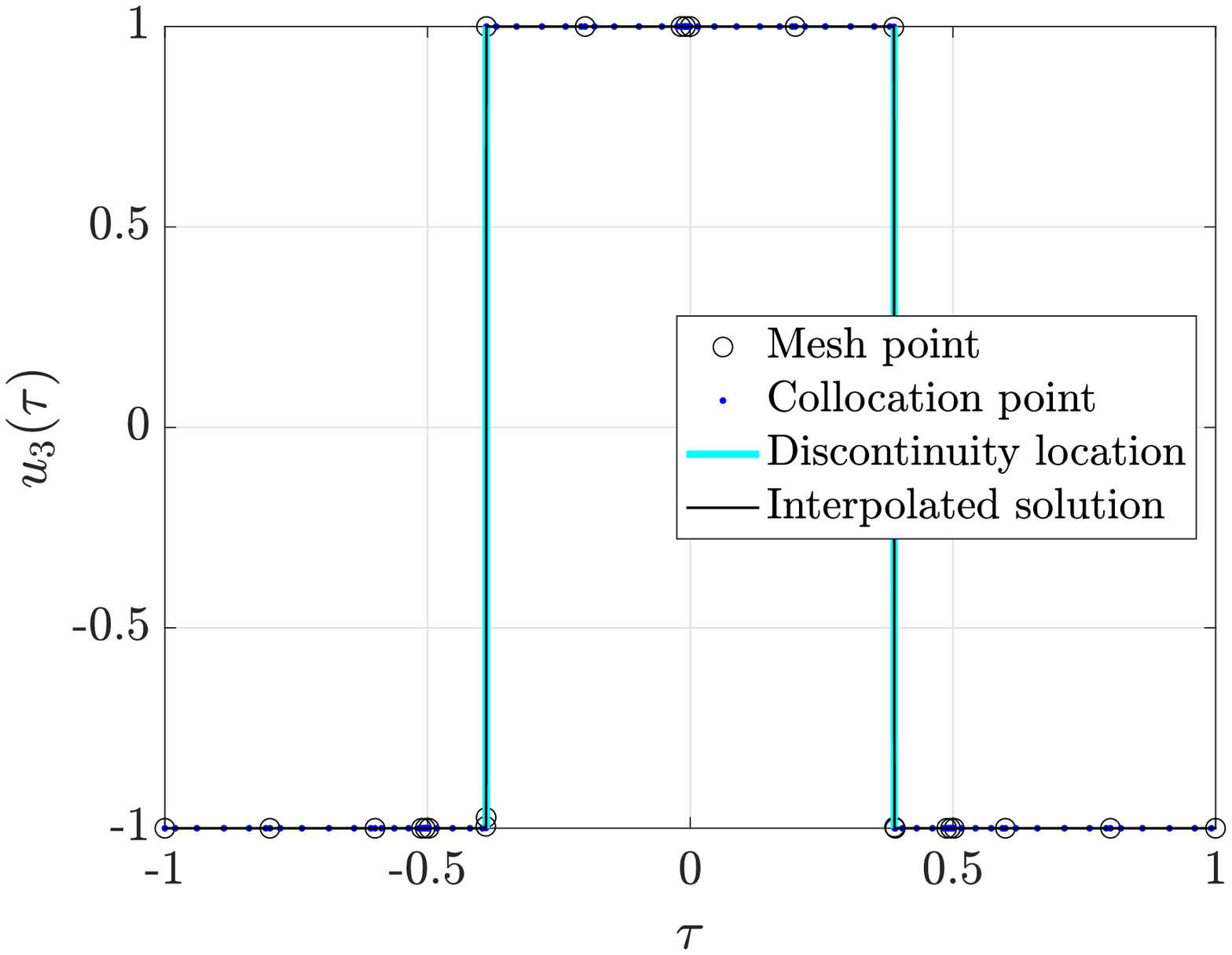}
  \caption{Control component solutions for Example 1 obtained by the $hp$-$(1)$ method at a mesh error tolerance of $\epsilon = 10^{-8}$.\label{fig:robotArm:control}}
\end{figure}

Consider the following optimal control problem obtained from Ref.~\citen{Munson} where the goal is to reorient a robotic arm in minimum time.  Minimize the final time
\begin{equation} \label{eq:robotArm-cost}
  \C{J} = t_f,
\end{equation}
subject to the dynamic constraints
\begin{equation} \label{eq:robotArm-dynamics}
\begin{array}{ll}
  \displaystyle \frac{d{y_1}}{d\tau} =  \frac{t_f - t_0}{2} y_2, &
  \displaystyle \frac{d{y_2}}{d\tau} =  \frac{t_f - t_0}{2 L}  u_1, \vspace{4pt} \\
  \displaystyle \frac{d{y_3}}{d\tau} =  \frac{t_f - t_0}{2}  y_4, &
  \displaystyle \frac{d{y_4}}{d\tau} =  \frac{t_f - t_0}{2 I_\theta} u_2, \vspace{4pt} \\
  \displaystyle \frac{d{y_5}}{d\tau} = \frac{t_f - t_0}{2}  y_6, & 	
  \displaystyle \frac{d{y_6}}{d\tau} =  \frac{t_f - t_0}{2 I_\phi}  u_3,
\end{array}
\end{equation}
the control inequality constraints
\begin{equation}\label{eq:robotArm-controlConstraint}
  -1 \leq u_i \leq 1 \quad (i = 1,2,3),
\end{equation}
and the boundary conditions
\begin{equation} \label{eq:robotArm-boundaryConditions}
\begin{array}{llcll}
  \displaystyle y_1(-1) = 9 / 2, & 	            \displaystyle y_1(+1) = 0, & \quad \quad & 
  \displaystyle y_2(-1) = 0, & 			        \displaystyle y_2(+1) = 0, \\
  \displaystyle y_3(-1) = 0, & 				    \displaystyle y_3(+1) = 2\pi / 3, & \quad \quad &
  \displaystyle y_4(-1) = 0, & 				    \displaystyle y_4(+1) = 0, \\
  \displaystyle y_5(-1) = \pi / 4, & 	        \displaystyle y_5(+1) = \pi / 4, & \quad \quad &
  \displaystyle y_6(-1) = 0, & 				    \displaystyle y_6(+1) = 0, \\
\end{array}
\end{equation} 
where $t_0 = 0$, $t_f$ is free, $I_\phi = \frac{1}{3} \left( (L-y_1)^3+y_1^3 \right)$, $I_\theta = I_\phi\sin^2(y_5)$, and $L = 5$.  A typical numerical solution for the control is shown in Fig.~\ref{fig:robotArm:control}.  As can be seen, the defining feature of the optimal control is its bang-bang structure with five discontinuities located at $\tau \approx \{-0.5000, -0.3882, 0.0000, 0.3882, 0.5000\}$.

\clearpage

\begin{figure}[h]
  \centering
  \begin{tabular}{cc}
  \includegraphics[height=2.25in]{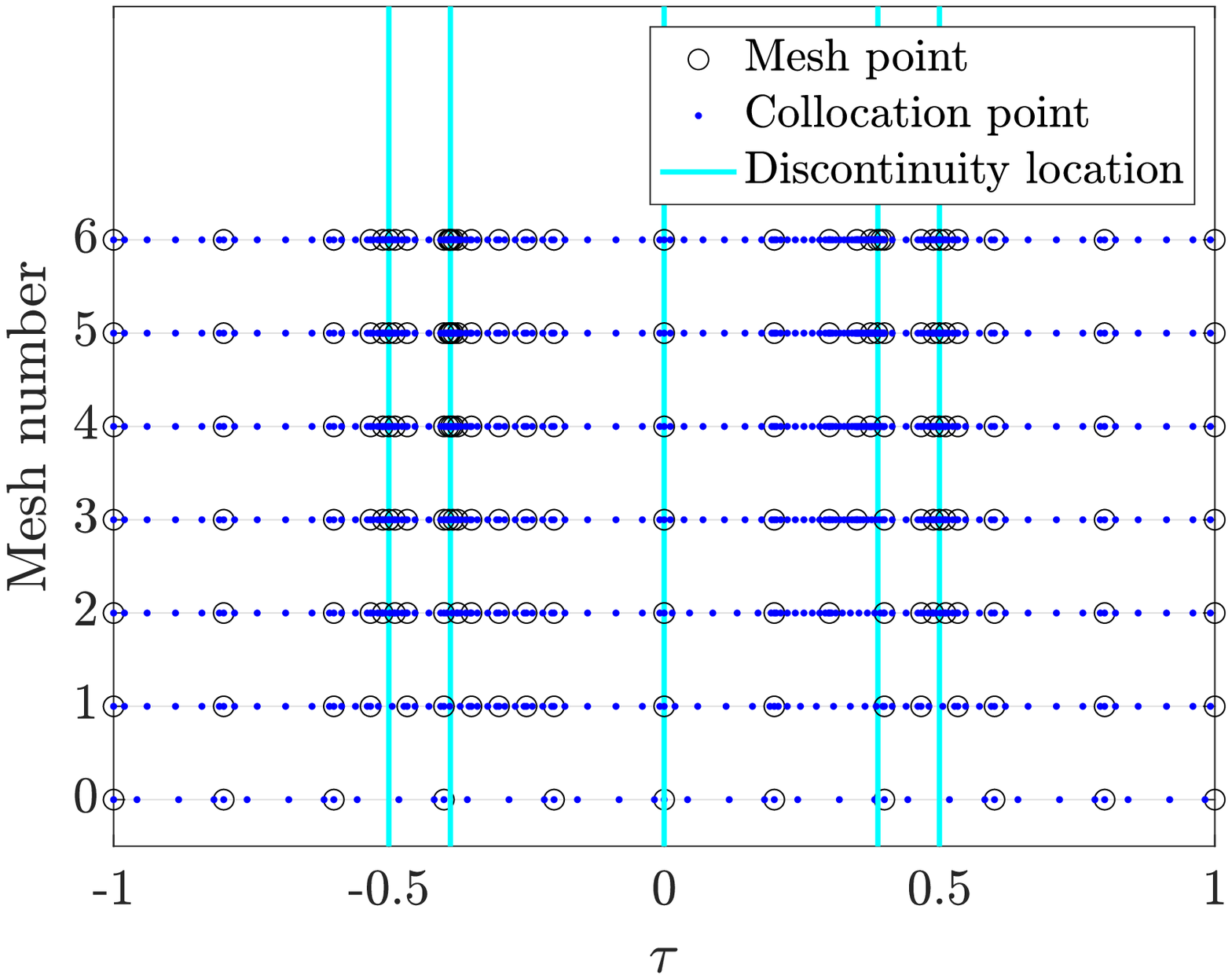}
  &
  \includegraphics[height=2.25in]{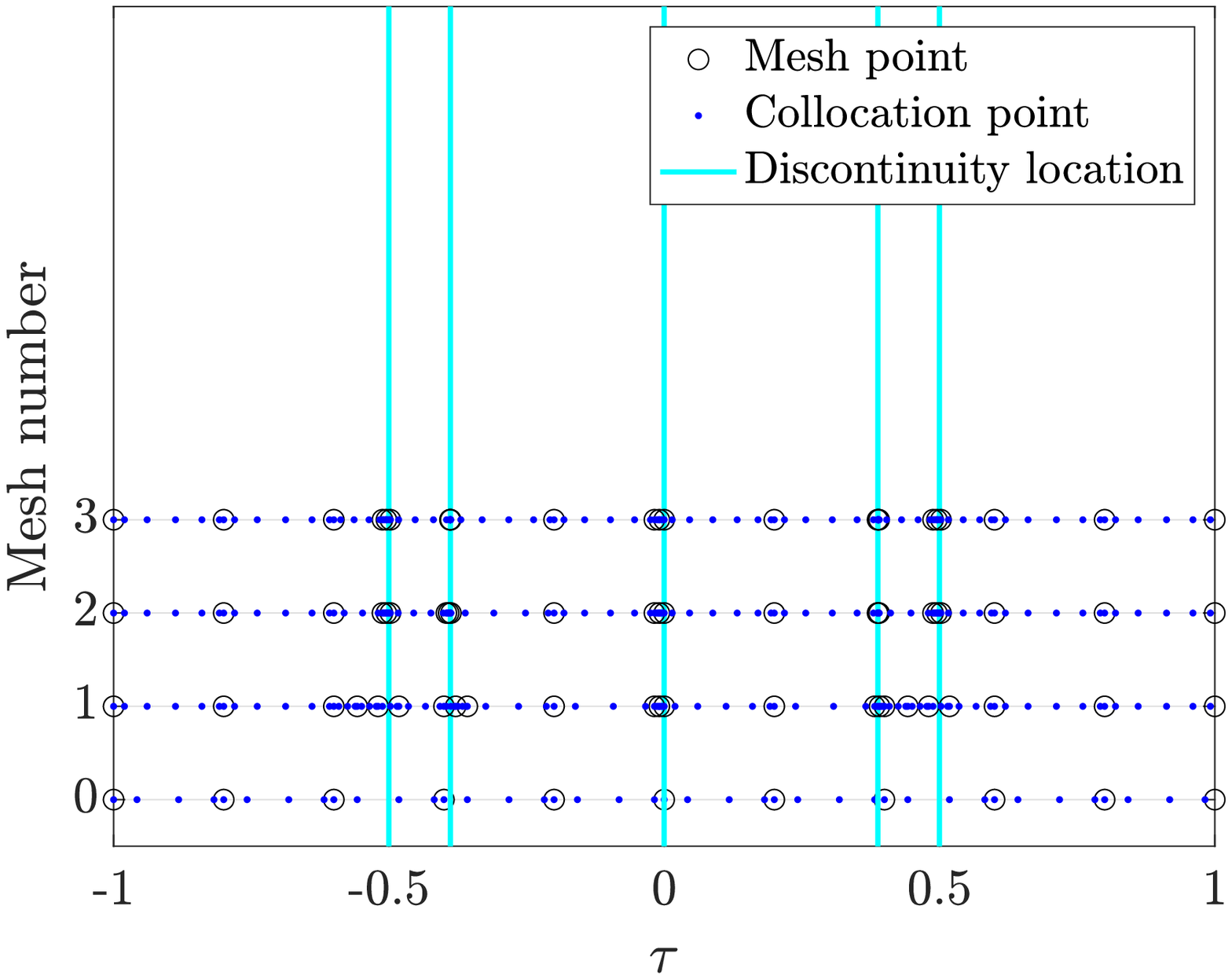}
  \end{tabular}
  \caption{Mesh histories obtained by the $hp$ method (left) and the $hp$-$(1)$ method (right) when solving Example 1 at a mesh error tolerance of $\epsilon = 10^{-8}$.\label{fig:robotArm:meshHistory}}
\end{figure}

The mesh histories obtained when solving this problem using the $hp$ and $hp$-$(1)$ methods are shown in Fig.~\ref{fig:robotArm:meshHistory} and are representative of other results obtained.  Observing the mesh history of the $hp$-$(1)$ method, it is seen that all five discontinuities were identified and bracketed by mesh points during the first mesh refinement iteration.  On subsequent refinement iterations, these brackets were updated to reflect the new discontinuity location estimates and their associated uncertainty bounds, each time producing mesh point brackets which bound the discontinuity locations in the optimal control solution.  The net effect is a systematic increase in resolution of the mesh around the discontinuity locations while keeping the mesh relatively sparse outside the neighborhoods of the discontinuities.  In contrast, the $hp$ method subdivides the mesh intervals containing discontinuities into equally spaced sub-intervals.  While such a strategy works in principle, it produces unnecessarily large meshes that concentrate grid points around the discontinuity locations.  In addition, the $hp$ method takes three more mesh refinement iterations than the $hp$-$(1)$ method to satisfy the desired solution accuracy.

\begin{figure}[H]
  \centering
  \begin{tabular}{rr}
  \includegraphics[height=2.25in]{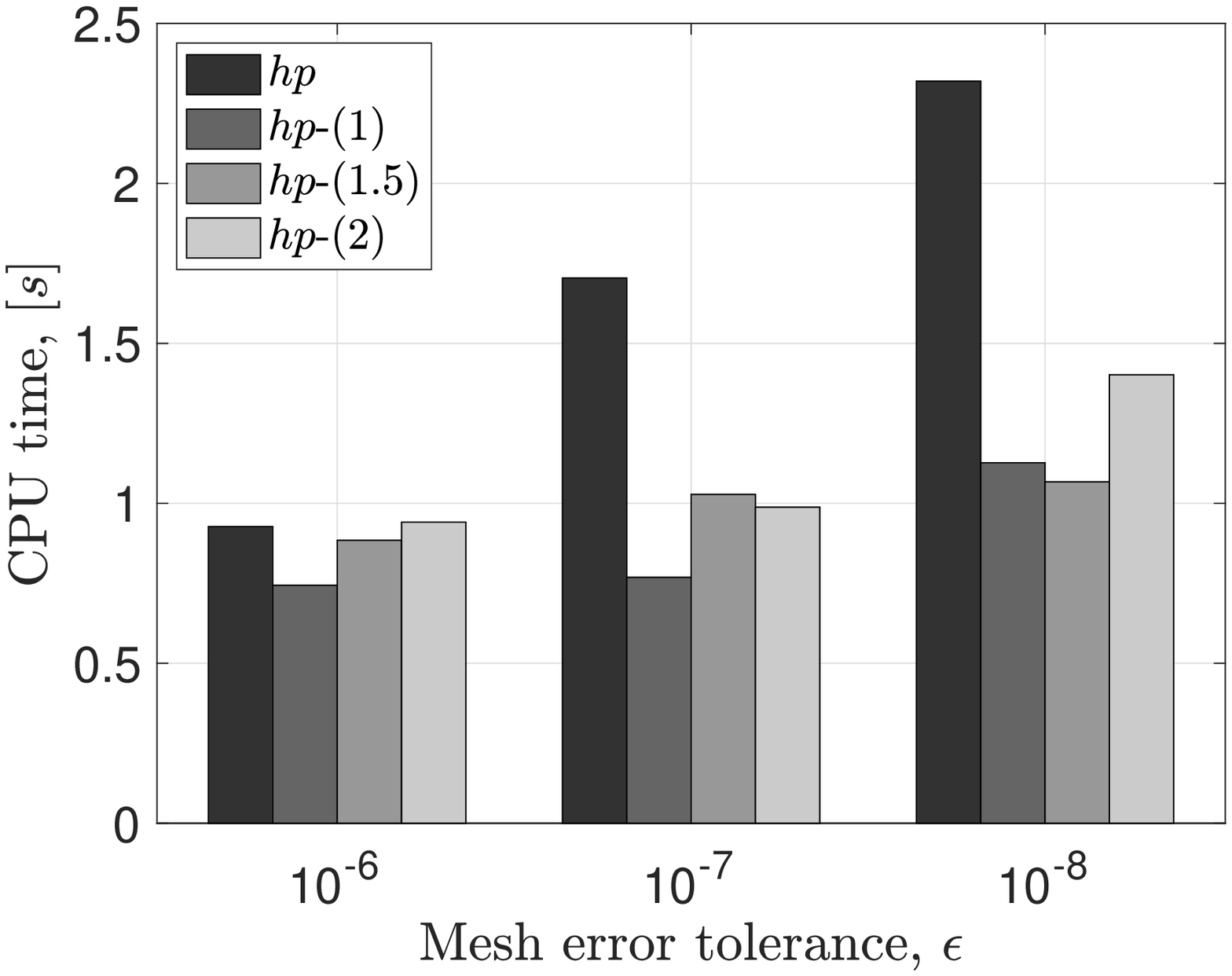}
  &
  \includegraphics[height=2.25in]{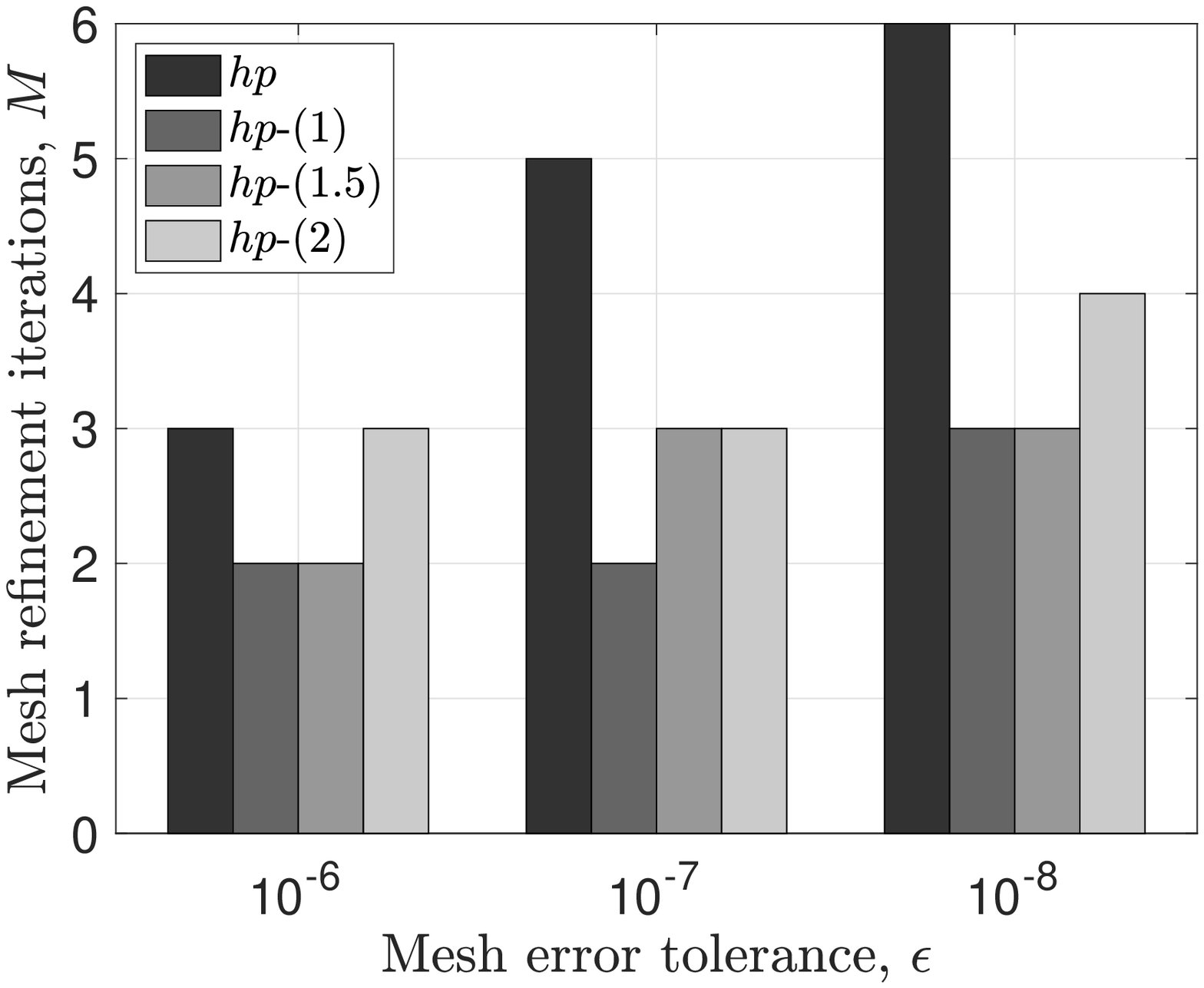}
  \vspace{1mm}
  \\
  \includegraphics[height=2.25in]{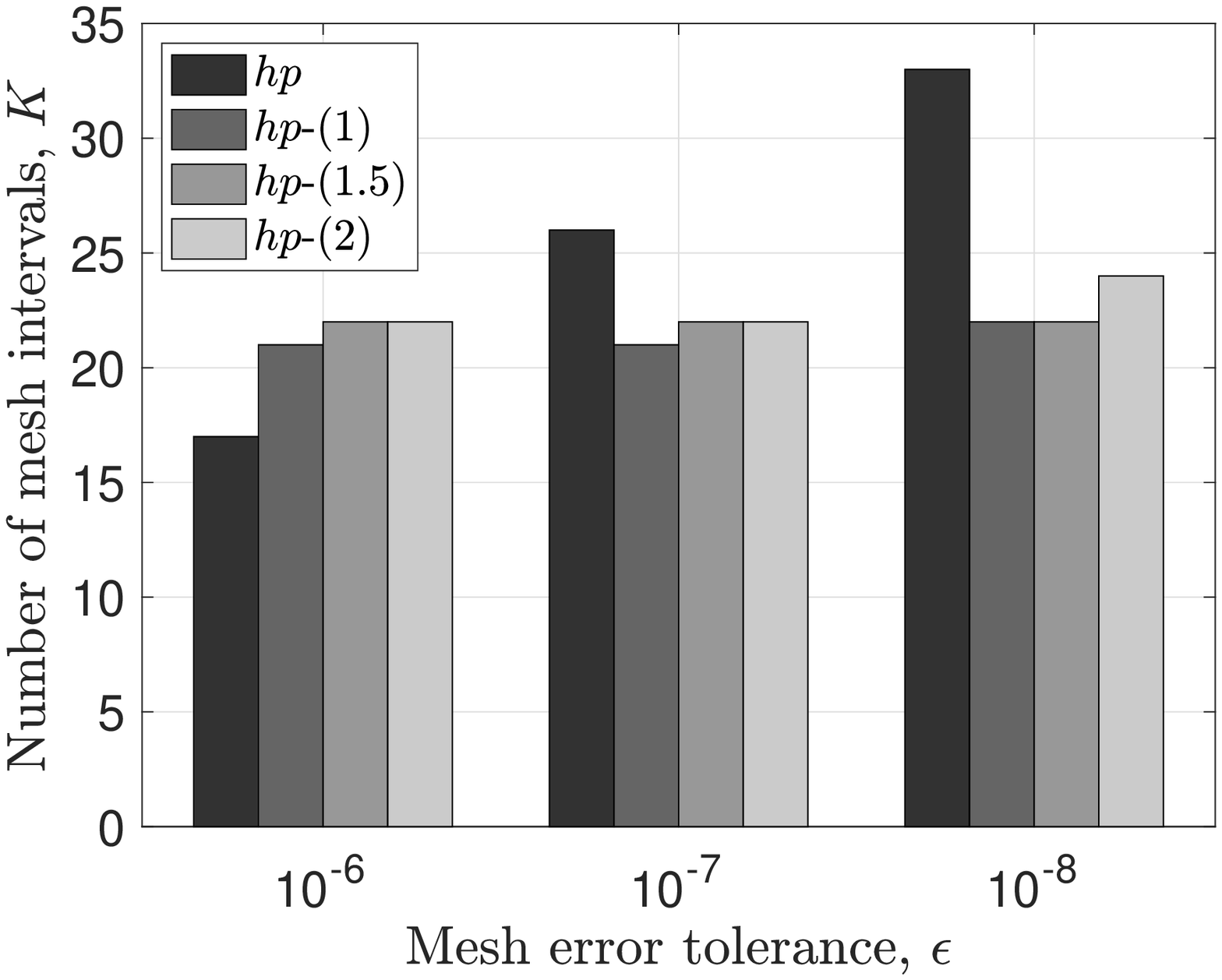}
  &
  \includegraphics[height=2.25in]{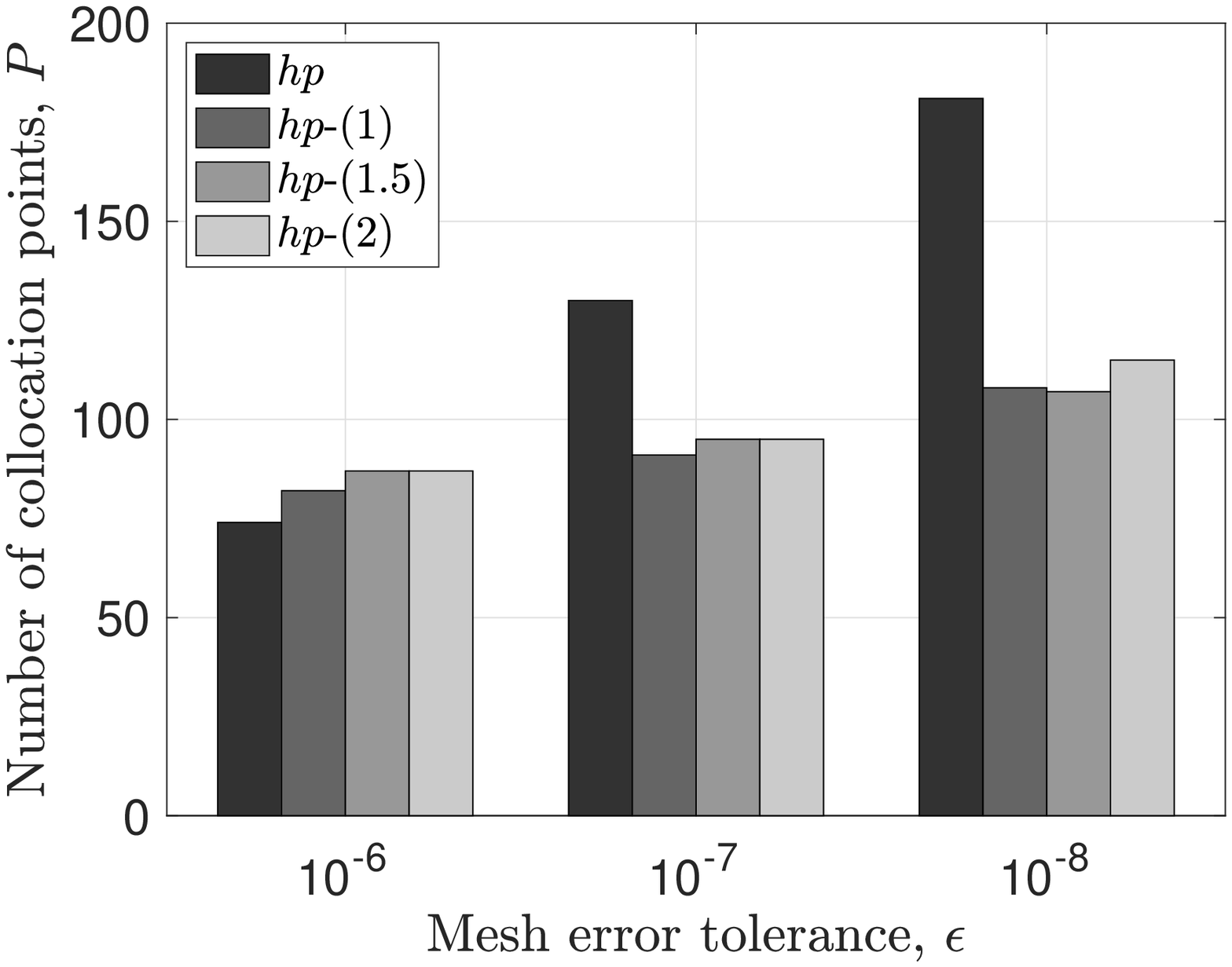}
  \end{tabular}
  \caption{Final mesh characteristics and computation times for Example 1.  \label{fig:robotArm:data}}
\end{figure}

Next, Fig.~\ref{fig:robotArm:data} shows the total computation time (CPU time) and the number of mesh refinement iterations ($M$) for each of the methods.  It is seen that the $hp$-$(1)$, $hp$-$(1.5)$, and $hp$-$(2)$ methods take the same or fewer number of mesh refinement iterations to converge to the solution accuracy when compared with the $hp$ method.  A similar trend is observed for the total computation time.  It is also observed that the largest differences in performance between the $hp$ method and the $hp$-$(1)$, $hp$-$(1.5)$, and $hp$-$(2)$ methods occurs at the two highest mesh error tolerances tested.  Comparing the $hp$-$(1)$, $hp$-$(1.5)$, and $hp$-$(2)$ methods among one another, it is observed that the number of mesh refinement iterations tends to increase slightly as the safety factor is increased from $\mu = 1$ to $\mu = 2$.  Figure~\ref{fig:robotArm:data} also shows the total number of mesh intervals ($K$) and the total number of collocation points ($P$) obtained on the final mesh for each of the methods.  Inspecting Fig.~\ref{fig:robotArm:data}, it is seen that the $hp$ method uses fewer collocation points and mesh intervals than the $hp$-$(1)$, $hp$-$(1.5)$, and $hp$-$(2)$ methods at the lowest error tolerance.  However, the opposite is true at the two highest mesh error tolerances.  In fact, the largest difference in performance occurs at the highest mesh error tolerance tested, where the $hp$-$(1)$, $hp$-$(1.5)$, and $hp$-$(2)$ methods use a much sparser mesh than the $hp$ method to achieve the same error tolerance of $\epsilon = 10^{-8}$.  Overall, Fig.~\ref{fig:robotArm:data} suggests that the nonsmooth mesh refinement method of Section~\ref{sect:meshRefinementMethod} is an effective approach for quickly and efficiently refining the mesh when discontinuities are present in the control solution, especially at high error tolerances.

\clearpage

% Reusable Launch Vehicle Entry-------------
\subsection{Example 2: Maximum Crossrange Shuttle Reentry\label{subsect:RLV}}
\begin{figure}[h]
  \centering
  \begin{tabular}{cc}
  \includegraphics[height=2.25in]{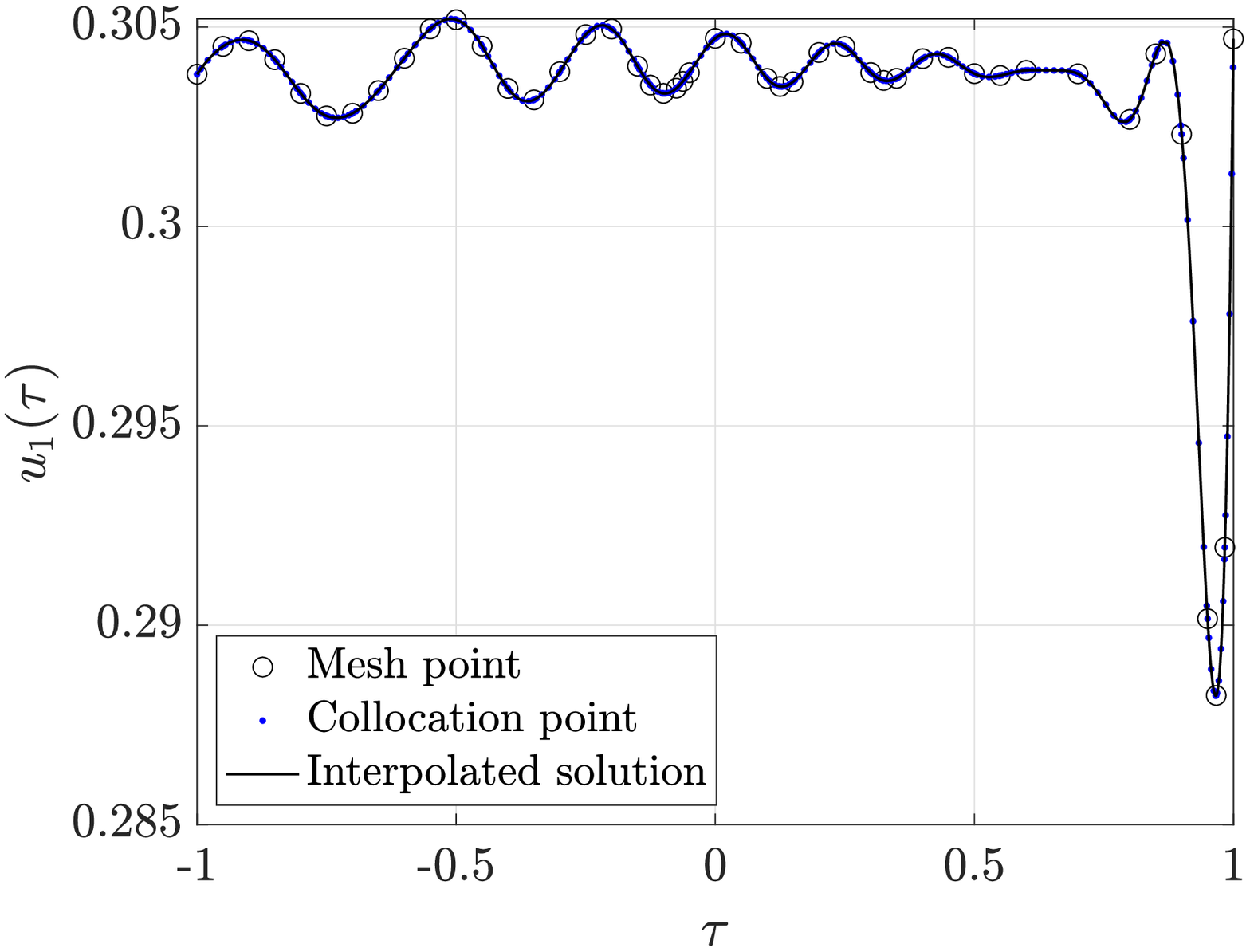}
  &
  \includegraphics[height=2.25in]{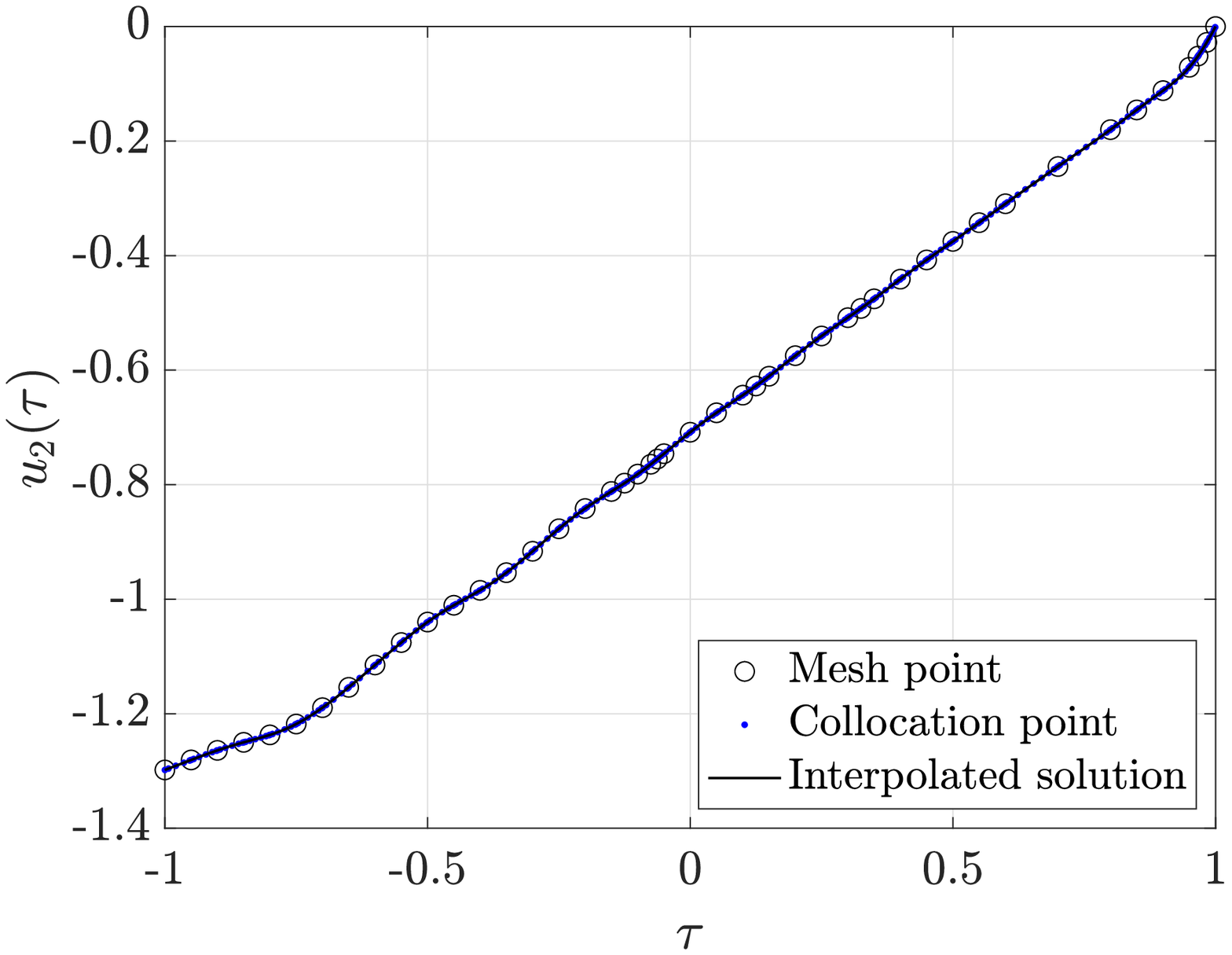}
  \end{tabular}
  \caption{Control component solutions for Example 2 obtained by the $hp$-$(1)$ method at a mesh error tolerance of $\epsilon = 10^{-8}$.\label{fig:rlv:control}}
\end{figure}

\begin{figure}[h]
  \centering
  \begin{tabular}{cc}
  \includegraphics[height=2.25in]{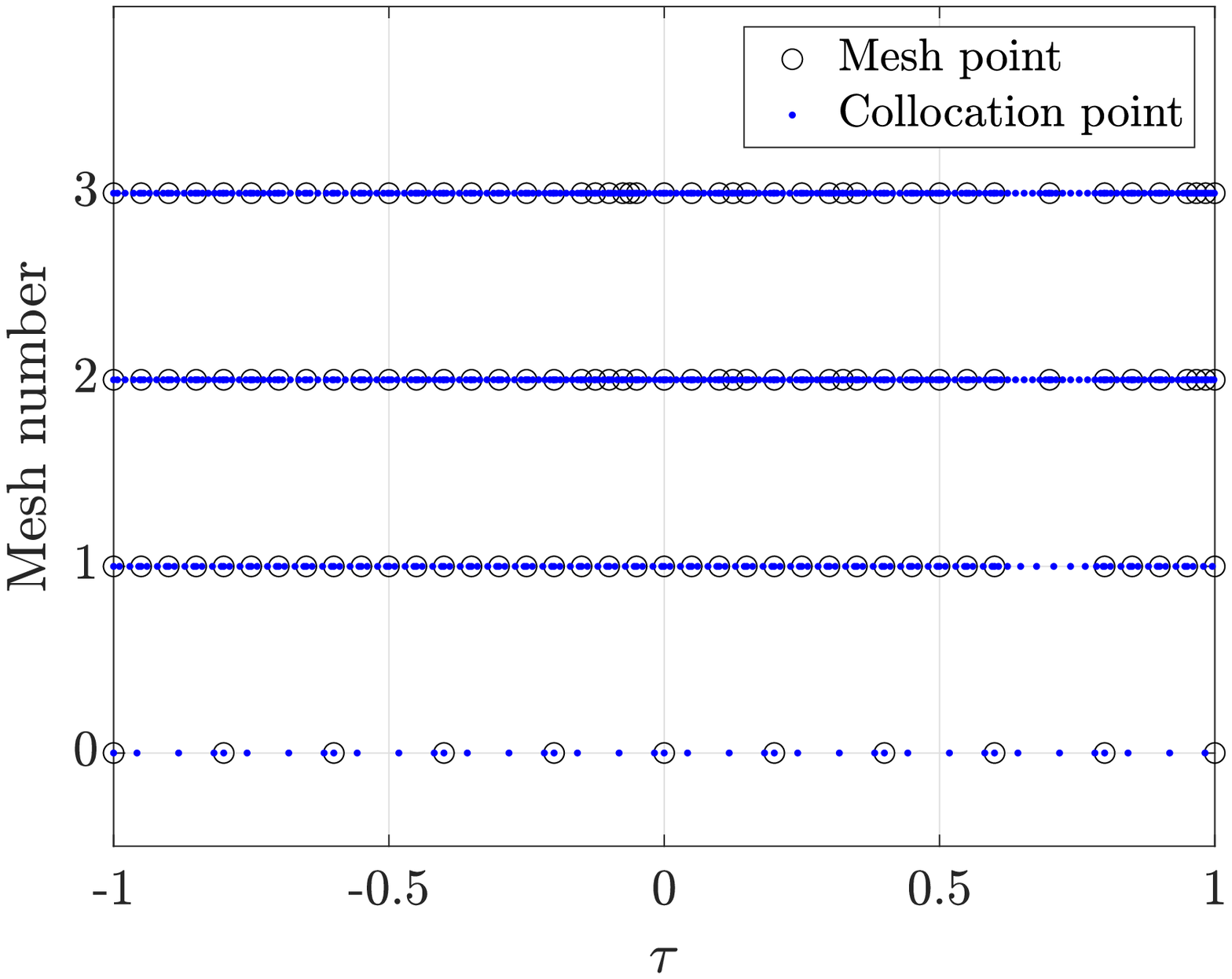}
  &
  \includegraphics[height=2.25in]{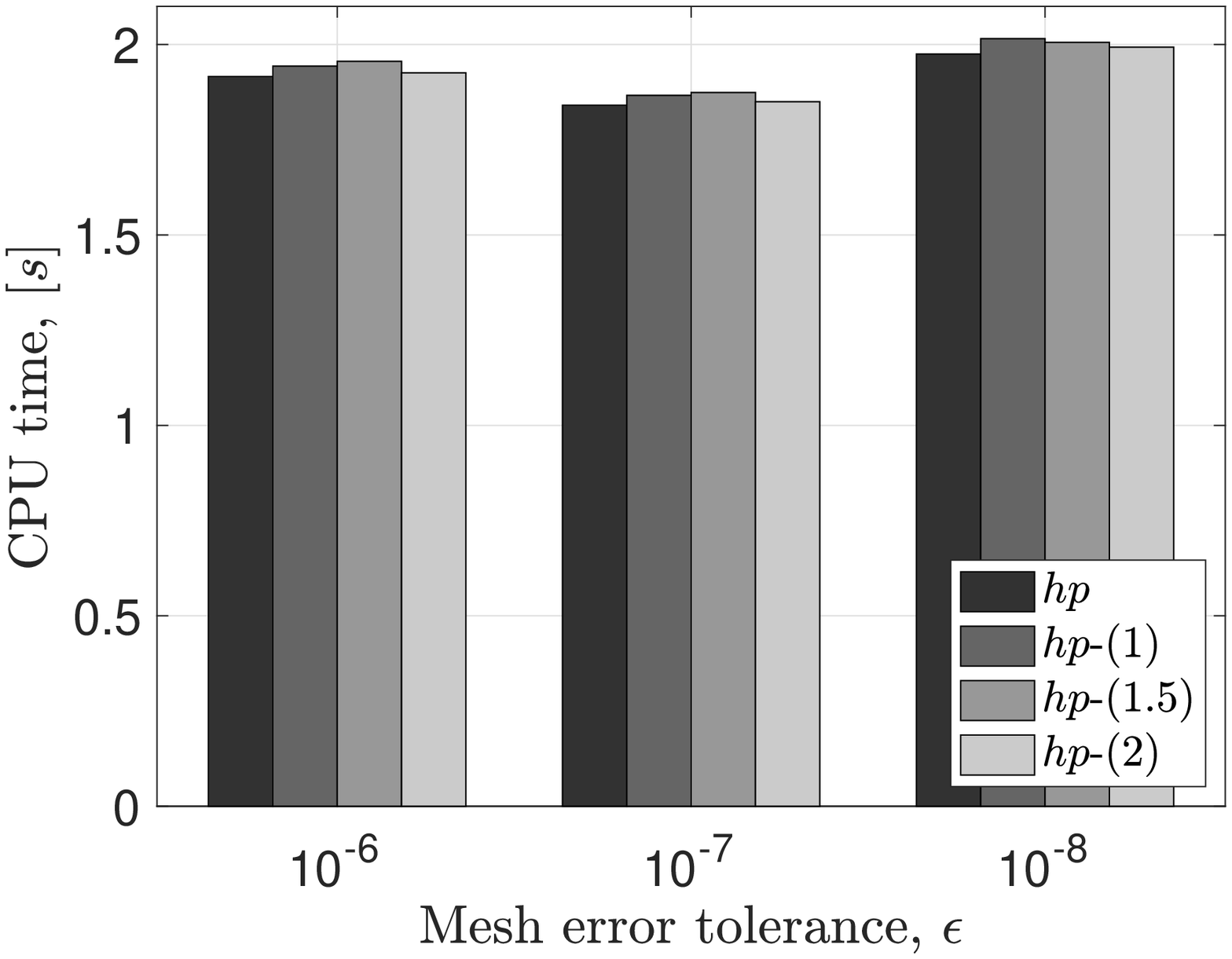}
  \end{tabular}
  \caption{Mesh history obtained by the $hp$, $hp$-$(1)$, $hp$-$(1.5)$, and $hp$-$(2)$ methods when solving Example 2 at a mesh error tolerance of $\epsilon = 10^{-8}$ (left).  Average computation times for Example 2 (right).  \label{fig:rlv:data}}
\end{figure}

Consider the reusable launch vehicle reentry problem of Ref.~\citen{Betts3}.  For brevity, the problem description is omitted from this paper, noting that Ref.~\citen{Betts3} uses English units whereas SI units were used to obtain the results shown here.  A typical numerical solution for the control is shown in Fig.~\ref{fig:rlv:control} where it is observed that the solution is continuous for both control components.  The continuous nature of the control solution allows one to observe how the nonsmooth mesh refinement method performs when no discontinuities are present in the control solution.  As one might hypothesize, the jump function approximations constructed for the two control components should indicate that no discontinuities are present and that only smooth mesh refinement is necessary.  If true, the mesh histories obtained using the $hp$, $hp$-$(1)$, $hp$-$(1.5)$, or $hp$-$(2)$ methods should be identical to one another with minimal deviations in the total computation times.  That hypothesis is confirmed for all three mesh error tolerances tested, and Fig.~\ref{fig:rlv:data} shows the identical mesh history obtained by all four methods as well as the average CPU times.  Thus, the $hp$-$(1)$, $hp$-$(1.5)$, and $hp$-$(2)$ methods correctly determine that the control solution is continuous, and they maintain the performance of the $hp$ method.

\vspace{2mm}

%----------------------------------------------
\section{Discussion\label{sect:discussion}}
Each of the examples in Section~\ref{sect:examples} highlights different features of the mesh refinement method developed in Section~\ref{sect:meshRefinementMethod}.  The first example demonstrates the ability of the method to accurately detect and locate discontinuities in the control solution as well as the method's ability to efficiently adjust the mesh around each identified discontinuity.  Moreover, the first example also demonstrates the computational benefits of using such an approach at high mesh error tolerances, reducing final mesh size and converging to the solution accuracy more quickly than the previously developed $hp$ method.  The second example highlights the ability of the method to discern when discontinuities are not present in the control solution.  In addition, it was seen that the $hp$-$(1)$, $hp$-$(1.5)$, and $hp$-$(2)$ methods revert back to the performance of the $hp$ method when no control discontinuities are present in the solution.  Thus, the method of Section~\ref{sect:meshRefinementMethod} is effective regardless of whether the solution is smooth or nonsmooth, because discontinuities are accurately identified and appropriate mesh refinement actions are taken in either case.

For problems with discontinuous control solutions, the $hp$-$(\mu)$ methods perform slightly differently depending on the value of the safety factor $\mu \geq 1$.  In the first example, it is seen that fewer mesh refinement iterations and shorter computation times can be achieved when a smaller value for $\mu$ is chosen.  Such a result makes sense, because the value of $\mu$ directly controls the widths of the two mesh intervals which bracket the discontinuity.  Decreasing $\mu$ produces narrower discontinuity brackets, increasing the resolution around the discontinuity more rapidly.  While a small discontinuity bracket is certainly desirable, it should be noted that smaller values of $\mu$ run a higher risk of generating discontinuity uncertainty bounds which miss the discontinuity, because the estimated bounds may be too narrow.

Two other parameters also affect the performance of the method of Section~\ref{sect:meshRefinementMethod}, but to a lesser degree.  The jump function approximation orders, $\mathscr{M} \subset \bb{N}^{+}$, and the discontinuity detection threshold, $0 < \eta < 1$, both affect \textit{whether} a discontinuity is detected.  Reference~\citen{Yoon2005} discusses choices for the jump function approximation orders, suggesting $\mathscr{M} = \{1,\ldots,6\}$ as a reasonable choice and noting that including $1 \in \mathscr{M}$ ensures first order convergence.  As for the discontinuity detection threshold, it determines the relative size of the jumps which should be detected.  The choice of $\eta$ is somewhat arbitrary but should be low enough to detect the discontinuities present in the solution yet high enough to filter out the continuous locations.  As the threshold is lowered, discontinuities are more likely to be identified.  Note, however, that regions of rapid change may be mistakenly identified as discontinuous if the threshold is lowered excessively.

Next, the mesh refinement method of Section~\ref{sect:meshRefinementMethod} is intended to be a step towards developing a general purpose $hp$-adaptive method capable of detecting, locating, and efficiently handling any type of discontinuity in the solution to an optimal control problem without apriori knowledge of the number of discontinuities or their locations.  In the ideal case, such a method would correctly identify all discontinuities of interest with no false positives and mitigate each discontinuity's detrimental effects on solution accuracy while keeping the mesh size and computation time as small as possible.  The method of Section~\ref{sect:meshRefinementMethod} could be further improved towards such a goal by employing other discontinuity detection techniques to locate nonsmooth solution behavior not captured by jump function approximations of the control components.  For example, jump function approximations built for the control derivative could help detect control derivative discontinuities.  Future research may also improve the accuracy and precision of the discontinuity location estimates on a given mesh, and could be more robust to noise in the solution.

%----------------------------------------------
\section{Conclusions\label{sect:conclusion}}
A mesh refinement method for solving optimal control problems has been developed.  The method has been shown to be capable of locating discontinuities in the control solution by employing jump function approximations.  These jump function approximations are generated using only the numerical solution for the control on a given mesh and require no apriori knowledge of the solution structure or discontinuity locations.  When discontinuities are identified, the mesh refinement method brackets the discontinuity locations with a specialized $h$-refinement approach.  The mesh refinement method has been demonstrated on two examples and compared against the previously developed method of Ref.~\citen{Liu2018}.  It was observed that the method required fewer mesh refinement iterations and achieved faster computation times when discontinuities were present in the control solution.

%----------------------------------------------
\section*{Acknowledgments}

The authors gratefully acknowledge support for this research from the U.S.~Office of Naval Research under grants N00014-15-1-2048 and N00014-19-1-2543, from the U.S.~National Science Foundation under grants DMS-1522629, DMS-1924762, and CMMI-1563225, and from the U.S.~Department of Defense (DoD) through the National Defense Science \& Engineering Graduate Fellowship (NDSEG) Program.

%----------------------------------------------
\renewcommand{\baselinestretch}{1}
\normalsize\normalfont

\bibliographystyle{aiaa}     % Number the references.
% \bibliography{references2}   % Use references.bib to resolve the labels.

\renewcommand{\baselinestretch}{1.5}
\normalsize\normalfont

\end{document}